\documentclass{amsart}

\usepackage{amsfonts,amssymb,amscd,amsmath,enumerate,verbatim, calc,xypic}

\CompileMatrices

\def\ds{\displaystyle}
\def\bu{\bullet}
\def\QED{\hfill\hbox{\qed}}

\def\ov{\overline}
\def\ul{\underline}

\def\wh{\widehat}
\def\col{\colon}
\def\dd{\partial}

\newcommand{\bsm}{{\boldsymbol m}}

\newcommand{\bsp}{{\boldsymbol p}}

\def\bC{{\boldsymbol C}}

\def\bG{{\boldsymbol G}}

\def\bK{{\boldsymbol K}}
\def\bL{{\boldsymbol L}}
\def\bM{{\boldsymbol M}}
\def\bN{{\boldsymbol N}}
\def\bP{{\boldsymbol P}}
\def\bQ{{\boldsymbol Q}}

\def\bS{{\boldsymbol S}}
\def\bT{{\boldsymbol T}}

\def\bZ{{\boldsymbol Z}}

\def\CC{{\mathcal C}}

\def\CGp{{\mathcal G\mathcal P}}
\def\CM{{\mathcal M}}

\newcommand{\BZ}{{\mathbb Z}}

\newcommand{\id}[1]{\operatorname {id}^{#1}}

\newcommand{\HH}[2]{\operatorname{H}_{#1}{#2}}
\newcommand{\Hh}[2]{\operatorname{H}^{#1}{#2}}

\newcommand{\Zh}[2]{\operatorname{Z}^{#1}{#2}}
\newcommand{\Ch}[2]{\operatorname{C}_{#1}{(#2)}}

\newcommand{\im}{\operatorname {Im}}
\newcommand{\maxx}{\operatorname{max}}
\newcommand{\Ker}{\operatorname{Ker}}
\newcommand{\Coker}{\operatorname{Coker}}

\newcommand{\Cone}{\operatorname{Cone}}

\newcommand{\xra}{\xrightarrow}
\newcommand{\xla}{\xleftarrow}
\newcommand{\lra}{\longrightarrow}

\newcommand{\depth}[2]{\operatorname{depth}_{#1}{#2}}

\newcommand{\idim}[2]{\operatorname{id}_{#1}{#2}}
\newcommand{\pd}[2]{\operatorname{pd}_{#1}{#2}}
\newcommand{\Gpd}[2]{\operatorname{Gpd}_{#1}{#2}}
\newcommand{\Gdim}[2]{\operatorname{G-dim}_{#1}{#2}}

\newcommand{\Hom}[3]{\operatorname{Hom}_{#1}({#2},{#3})}

\newcommand{\Ext}[4]{\operatorname{Ext}_{#1}^{#2}(#3,#4)}
\newcommand{\Extt}[4]{\operatorname{E\overset{\scriptscriptstyle\wedge}{x}t}\!{\vphantom
E}_{#1}^{#2}(#3,#4)}

\newcommand{\Extr}[3]{\operatorname{Ext}_{\CGp}^{#1}(#2,#3)}

\newcommand{\etH}[4]{\operatorname\eth\!^{#1}_{#2}(#3,#4)}
\newcommand{\etr}[3]{\operatorname{\eth}\!_{\CGp}^{#1}(#2,#3)}
\newcommand{\ett}[4]{\operatorname{\overset{\scriptscriptstyle\wedge}\eth}\!^{#1}_{#2}(#3,#4)}

\newcommand{\epsr}[3]{\operatorname{\varepsilon}\!_{\CGp}^{#1}(#2,#3)}
\newcommand{\epst}[4]{\operatorname{\overset{\scriptscriptstyle\wedge}\varepsilon}\!_{#1}^{#2}(#3,#4)}

\newcommand{\deltat}[4]{\operatorname{\overset{\scriptscriptstyle\wedge}\delta}\!_{#1}^{#2}(#3,#4)}

\theoremstyle{plain}
\newtheorem{theorem}{Theorem}[subsection]

\newtheorem{lemma}[theorem]{Lemma}
\newtheorem{proposition}[theorem]{Proposition}

\newtheorem{5theorem}{Theorem}[section]
\newtheorem{5proposition}[5theorem]{Proposition}
\newtheorem{5corollary}[5theorem]{Corollary}
\newtheorem{5lemma}[5theorem]{Lemma}
\newtheorem{5propdef}[5theorem]{Proposition-Definition}

\theoremstyle{definition}

\newtheorem{chunk}[theorem]{}

\newtheorem{5chunk}[5theorem]{}
\newtheorem{5example}[5theorem]{Example}
\newtheorem{5remark}[5theorem]{Remark}

\newtheorem{5definition}[5theorem]{Definition}
\newtheorem{5construction}[5theorem]{Construction}

\theoremstyle{remark}

\newtheorem{step}{Step}

\begin{document}

\date{\today}

\title{Gorenstein projective dimension for complexes}

\author[O. Veliche]{Oana Veliche}
\address{Department of Mathematics, Purdue University, West Lafayette,
Indiana~47907}
\email{oveliche@math.purdue.edu}

\begin{abstract} 
We define and study a notion of Gorenstein projective dimension for
complexes of left modules over associative rings. For complexes of
finite Gorenstein projective dimension we define and study a Tate
cohomology theory. Tate cohomology groups have a natural transformation
to classical Ext groups.  In the case of module arguments, we show that
these maps fit into a long exact sequence, where every third term is a
relative cohomology group defined for left modules of finite Gorenstein
projective dimension.
 \end{abstract}

\maketitle

\section*{Introduction}
  
We study generalized homological dimensions for complexes
of modules. More precisely we define and study a notion of Gorenstein
projective dimension for complexes of left modules over associative rings.

In the classical book of Cartan and Eilenberg \cite{ce} concepts of
projective, injective and weak (flat) dimensions were defined for left
modules over arbitrary rings. New dimensions have been defined since then.
For finite modules over commutative noetherian rings Auslander and Bridger
\cite{ab} introduced a Gorenstein dimension.  The reason for the name is
that a commutative local ring  is Gorenstein if and only if every finite
module has finite Gorenstein dimension.  More recently, several dimensions
have been defined over commutative noetherian local rings: complete intersection
dimension by Avramov, Gasharov and Peeva \cite{avgp}, polynomial complete
intersection dimension by Gerko \cite{g}, upper Gorenstein projective
dimension  by Veliche \cite{v} and  Cohen-Macaulay dimension by Gerko
\cite{g}. These dimensions are related to the ring in the same way as
the Gorenstein dimension is related to the Gorenstein rings.  An overview
of these dimensions can be found in \cite{av}.

For every left module $M$ over an associative ring $R$, Enochs and Jenda
\cite{ej1} defined a Gorenstein projective dimension denoted $\Gpd RM$;
they studied it when the ring is  coherent or $n$-Gorenstein. For finite
modules over commutative noetherian rings it coincides with the Auslander
and Bridger's Gorenstein dimension. For left modules over arbitrary
associative rings, Holm \cite{hh} proved that the new concept has the
desired properties.

In a different direction, homological dimensions have been extended to
complexes. Avramov and Foxby \cite{af} defined projective, injective, and
flat dimensions for arbitrary complexes of left modules over associative
rings.  Over commutative local rings, Yassemi \cite{y} and Christensen
\cite{cl} introduced a Gorenstein projective dimension for complexes with
bounded below homology.  Complete intersection dimension was extended
to homologically finite complexes by Sather-Wagstaff \cite{sw}.

The main purpose of this paper is to introduce and study a concept of
Gorenstein projective dimension $\Gpd R\bM$ associated to every (not necessarily
bounded) complex $\bM$ of left modules over an arbitrary associative
ring $R$.  It is obtained by blending the approach of Enochs and Jenda
to Gorenstein projective dimension of modules with that of Avramov and
Martsinkovsky \cite{am} to Gorenstein dimension of finite modules over noetherian rings.
We show that the  previously defined notions, when they can be applied,
agree with ours, and that most properties  of modules of finite
Gorenstein projective dimension are preserved for complexes.  Next
we describe our construction and some results in more detail.

A complex $\bT$ is called \emph{totally acyclic} if its modules
are projective, it is exact, and the complex  $\Hom{R}{\bT}{Q}$
is exact for every projective module $Q$.  We say that a complex
$\bM$ has $\Gpd R\bM\le g$ if and only if there exist morphisms of complexes $\bT\xra{\tau}\bP\xra{\pi}\bM$, where is $\bT$ totally acyclic, $\tau_i$ bijective for all $i\geq g$ (with the grading convention from \ref{complexes}), and $\bP\xra{\pi}\bM$ is a semiprojective resolution (in the sense of
\cite{afh}); such a diagram is called a \emph{complete resolution} of $\bM$.
In \ref{complete} we prove that the Gorenstein projective dimension of
$\bM$ can be computed from  any semiprojective resolution $\bP\xra{}\bM$.
In \ref{Gpd_pd} we relate it to projective dimension by an inequality
$\Gpd R\bM\leq\pd R\bM$, where equality holds when the right hand side
is finite.  When $\Gpd R\bM$ is finite, we prove in \ref{complete_Gpd}
that it is equal to the least integer $n$ for which $\Ext{R}{i}{\bM}{Q}=0$
for all $i\ge n$ and all projective modules $Q$.

Tate cohomology for finite modules of finite Gorenstein dimension
over noetherian rings is treated, most recently, by Avramov and
Martsinkovsky in \cite{am}. In Section \ref{Tate cohomology} we extend their construction to
the case of complexes and show that most of its properties for modules
extend well to the framework of complexes.  For each complex $\bM$
of finite Gorenstein projective dimension and every complex $\bN$, the Tate cohomology $\Extt{R}{}{\bM}{\bN}$ is the cohomology
of the complex $\Hom{R}{\bT}{\bN}$, where $\bT\xra{}\bP\xra{}\bM$ is
a complete resolution. It is easily seen to be a covariant functor of
the second argument.  In \ref{tate_contravariant} we prove that it is a
contravariant functor in the first argument; this is rather more delicate.
As in the case of modules, Tate cohomology is rigid: In \ref{tate_pd} we
prove that $\Extt{R}{i}{\bM}{-}=0$ for a single $i\in\BZ$ is equivalent
to $\pd R\bM<\infty$.  If $\bM$ is bounded below, then we also prove
that $\pd R\bM< \infty$ if and only if $\Extt{R}{i}{-}{\bM}=$ for some
(any) $i\in\BZ$, if and only if $\Extt{R}{0}{\bM}{\bM}=0$.

In Section \ref{Complexes with bounded below homology}, we show that
when the complex $\bM$ has bounded homology, $\Gpd R\bM$ can be
expressed in terms of finite and of special Goresntein projective
resolutions, see \ref{Gpd_finite_special}. In particular, we show in
\ref{old_new} that our definition of Gorenstein
projective dimension agrees with the ones of Gorenstein (projective) dimension of
Christensen and Yassemi when they can be applied.   

In the last section, we show the functorial properties of relative
cohomology for left modules over associative rings and show in
\ref{relative_Tate} that for any module of finite Gorenstein projective
dimension there exists a long exact sequence relating the natural Ext,
Tate cohomology and relative cohomology.

\section{Complexes}
\label{Complexes}

In this section we recall basic definitions and properties of complexes used in this paper.
For more details the reader can consult \cite{afh} or \cite{iv}. 

\subsection{Complexes}
\label{complexes}

Let $R$ be an associative ring, $\CM=\CM(R)$ the category of left
$R$-modules  and $\CC=\CC(R)$ the category of complexes of left
$R$-modules.

To every complex
$$
\bM=\quad \cdots\rightarrow
M_n\xrightarrow{\partial^\bM_n}M_{n-1}\rightarrow\cdots
$$
in $\CC$ we associate the numbers  
$$
\sup\bM=\sup\{i\mid M_i\not=0\}\quad\text{and}\quad
\inf\bM=\inf\{i\mid M_i\not=0\}.
$$
The complex $\bM$ is called \emph{bounded above} when $\sup\bM<\infty$, \emph{bounded below} when $\inf\bM>-\infty$ and \emph{bounded} when it is bounded below and above. 

The $n^{\text{th}}$ \emph{homology module} of $\bM$ is the module
$\HH{n}{(\bM)}=\Ker(\partial^\bM_n)/\ds\im(\partial^\bM_{n+1})$; we set $\Hh{n}{(\bM)}=\HH{-n}{(\bM)}$.
We also set $\Ch{n}{\bM}=\Coker(\partial_{n+1}^\bM)$.  The $i^{\text {th}}$ \emph{shift} of $\bM$ is a complex
$\Sigma^i\bM$ with $n^{\text{th}}$ component $M_{n-i}$ and differential
$\partial_n^{\Sigma^i\bM}=(-1)^i\partial^\bM_{n-i}$; we write $\Sigma \bM$ instead of $\Sigma^1 \bM$.

In the following discussion, $\bM$ and $\bN$ denote complexes of left $R$-modules.

\begin{chunk}   
A homomorphism $\varphi\col \bM\to\bN$ of
degree $n$ is a family $(\varphi_i)_{i\in\BZ}$ of homomorphisms of $R$-modules
$\varphi_i\col M_i\to N_{i+n}$.  All such homomorphisms form an
abelian group, denoted $\Hom{R}{\bM}{\bN}_n$; it is clearly isomorphic
to $\prod_{i\in\BZ}\Hom{R}{M_i}{N_{i+n}}$.  We let $\Hom{R}{\bM}{\bN}$
denote the complex of $\BZ$-modules with $n^{\text{th}}$ component
$\Hom{R}{\bM}{\bN}_n$ and differential
$$
\partial(\varphi)=\partial ^\bN \varphi-(-1)^{|\varphi|}\varphi\partial^\bM.
$$
A homomorphism $\varphi\in\Hom{R}{\bM}{\bN}_n$ is called a
\emph{chain map} if $\partial(\varphi)=0$, i.e., if 
$$
\partial ^\bN_{i+n} \varphi_i=(-1)^{|\varphi|}\varphi_{i-1}\partial^\bM_i
\quad\text{for all}\quad i\in\BZ.
$$
A chain map of degree 0 is called \emph{morphism}.
Homomorphisms $\varphi$ and $\varphi'$ in $\Hom{R}{\bM}{\bN}_n$ are called
\emph{homotopic}, denoted $\varphi \sim \varphi'$, if there exists a degree
$n+1$ homomorphism $\varkappa$, called a \emph{homotopy}, such that
$\partial(\varkappa)=\varphi-\varphi'$.\  A \emph{homotopy equivalence} is a morphism $\varphi\col\bM\to\bN$ for which  there exists a morphism $\psi\col\bN \to\bM$ such that $\varphi\psi \sim \id{\bN}$ and $ \psi\varphi\sim \id{\bM}$. 
\end{chunk}

\begin{chunk}
\label{hom_spliting}
\cite[(6.2.7)]{afh} Every morphism $\varphi\col\bM\to\bN$ has a factorization  
$$
\bM\xra{\mu}\bM'\xra{\varphi'}\bN
$$
where 
$
\bM'=\bM\oplus\Cone (\id{\bN}),
$
$\mu$ is a homotopy equivalence and $\varphi '$ is a surjective morphism. 
\end{chunk}

\begin{chunk}  
\label{contractible}
The complex $\bM$ is called \emph{contractible} if it is  homotopy  equivalent to ${\bf 0}$. 
This is the case if and only if $\Hom{R}{\bM}{\bM}$ is exact, if and only if $\bM\cong\Cone(\id \bZ)$ where $\bZ$ is a complex with trivial differential, see \cite[(6.2.8)]{afh}. 
\end{chunk}

\begin{chunk} 
A \emph{quasiisomorphism} $\varphi\col \bM\to \bN$ is a morphism, such that the
induced map $\HH{n}(\varphi)\col\HH{n}(\bM)\to\HH{n}(\bN)$ is an isomorphism
for all $n\in\BZ$. The complexes $\bM$ and $\bN$ are called
\emph{quasiisomorphic}, and denoted by $\bM\simeq\bN$, if they can be linked by a sequence of quasiisomorphisms with arrows in alternating directions.
\end{chunk}

\begin{chunk}
\cite[(1.1.1)]{am}
\label{quism}
 Let $\beta\col \bN\to\bN'$ be a
morphism of complexes such that $\Hom{R}{\bM}{\beta}\col
\Hom{R}{\bM}{\bN}\to\Hom{R}{\bM}{\bN'}$ is a quasiisomorphism.  For each
morphism $\gamma\col\bM\to \bN'$ there is a morphism
$\alpha\col \bM\to\bN$ with $\gamma\sim \beta\,\alpha$ (even
$\gamma=\beta\,\alpha$ if $\Hom{R}{\bM}{\beta}$ is surjective).  If
$\gamma'\col \bM\to \bN'$ and $\alpha'\col \bM\to\bN$ satisfy
$\gamma'\sim\gamma$ and $\gamma'\sim \beta\,\alpha'$, then
$\alpha'\sim\alpha$.
\end{chunk}

From \cite[(\S 5, no. 2, Proposition 2) ]{b} and \cite[(4.5.3.1), (5.2.3.1)]{afh} we obtain.

\begin{chunk}
\label{proj_exact}
 If $\bP$ is a complex of projectives, then $\Hom{R}{\bP}{-}$ preserves  quasiisomorphisms between bounded below complexes.
  \end{chunk}

\subsection{Connecting homomorphisms}
\label{theta_maps}
Let 
$
\bM_\bu=\  0\to \bM\xra{\mu}\bM'\xra{\mu'}\bM''\to 0
$
be an exact sequence of complexes which splits as  sequence of graded modules. There
exist homomorphisms $\nu''\col\bM''\to\bM'$ and $\nu'\col\bM'\to\bM$ such that 
$$ 
\nu'\nu''=0,\quad \mu'\nu''=\id{\bM''},\quad
\nu'\mu=\id{\bM}\quad\text{and}\quad \nu''\mu'+\mu\nu'=\id{\bM'}.
$$

Set $\theta^{\bM''}=\nu'\dd^{\bM'}\nu''$.

\begin{chunk}
\label{theta_morphism}
\cite[(4.8), p.19]{iv} 
The map $\theta^{\bM''}\col\bM''\to\Sigma\bM$ is a morphism of complexes and
$\eth=\HH{\!}{(\theta^{\bM''})}$, where $\eth$ is the connecting homomorphism
of the cohomology long exact sequence of $\bM_\bu$.
\end{chunk}

\begin{chunk}
\cite[(4.9), p.20]{iv}
\label{theta_arg2}
For any complex $\bN$ the connecting homomorphism of the long exact homology sequence
$$
0\to\Hom R{\bM''}{\bN}\to\Hom{R}{\bM'}{\bN}\to\Hom{R}{\bM}{\bN}\to 0
$$
is given by 
$$
\etH{n-1}{R}{\bM_\bu}{\bN}=
(-1)^{n}\Hh{n}{(\Hom{R}{\theta^{\bM''}}{\bN})}
\quad \text{for all}\quad n\in\BZ.
$$
\end{chunk}

\begin{chunk}
\label{theta_construction}
Let $\lambda\col\bM\to\bN$ be a morphism of complexes.  From the
exact sequence 
$$
0\to\bN\xra{}\Cone{(\lambda)}\xra{}\Sigma \bM\to 0
$$ 
we get a short exact sequence which splits
as a sequence of graded modules
$$
0\to\Sigma^{-1}\bN\xra{\mu}\Sigma^{-1}\Cone{(\lambda)}\xra{\mu'}\bM\to 0
$$
Let $\nu''\col\bM\to\Sigma^{-1}\Cone{(\lambda)}$ be the canonical
injection and $\nu'\col\Sigma^{-1}\Cone{(\lambda)}\to\Sigma^{-1}\bN$
be the canonical surjection. Then, the map $\theta^\bM$  defined in
\ref{theta_maps} is given by
$$
\theta^\bM=\lambda.
$$
\end{chunk}

\subsection{Semiprojective resolutions} 

Let $\bP$ be a complex of $R$-modules.

Following \cite{afh}, we say that $\bP$ is {\it
semiprojective} if $\Hom{R}{\bP}{-}$ preserves surjective
quasiisomorphisms.

\begin{chunk}
\cite[(9.5.1), (9.4.1)]{afh}
\label{semiprojective} 
The following conditions are equivalent.
\begin{enumerate}[\rm\quad(i)]
\item $\bP$ is semiprojective.
\item Each $P_i$ is projective and $\Hom{R}{\bP}{-}$
preserves quasiisomorphisms.
\item Each $P_i$ is projective and $\HH{\!}{(\Hom{R}{\bP}{\bN})}=0$
for every complex $\bN$ with $\HH{\!}{(\bN)}=0$.
\end{enumerate}
\end{chunk}

\begin{chunk}
\cite[(1.4.P)]{af}
\label{semiprojective_quism} 
If $\bP$ is a semiprojective complex and $\bP\simeq\bM$, then there
exists a quasiisomorphism $\bP\to\bM$.
\end{chunk}

From \ref{contractible} and \ref{semiprojective} we get the following.

\begin{chunk} 
\label{semiproj_contr}
If $\bP$ is a semiprojective complex with $\HH{}{(\bP)}=0$, then $\bP$
is contractible.
\end{chunk}

A \emph{semiprojective resolution of $\bM$} is a quasiisomorphism of
complexes $\pi\col \bP\to\bM$, with $\bP$ semiprojective; when $\pi$
is surjective, the resolution is called \emph{strict}.

\begin{chunk}
\cite[(9.3.2), (12.3.2), (12.2.7)]{afh}
\label{semi_res}
Every complex $\bM$ has a strict semiprojective resolution
$\bP\to\bM$.  If $\HH{\!}{(\bM)}$ is bounded below, then $\bP$
can be chosen so that $\inf\bP=\inf\HH{\!}{(\bM)}$. If, in addition,
$R$ is left noetherian and $\HH i{(\bM)}$ is finitely generated for
each $i\in\BZ$, then $\bP$ can be chosen with each $P_i$ finitely
generated.
\end{chunk}

\begin{chunk}
\cite[(9.7.3.2')]{afh}
\label{proj}
If $M$ is a left module, then  any classical projective resolution
$\bP\to M$ is a strict semiprojective resolution.
\end{chunk}

The next two results were communicated to the author by Luchezar
Avramov. They are special cases of theorems that will eventually be
included in \cite{afh}.

\begin{proposition}[Shanuel's Lemma for Complexes]
\label{Shanuel}
If $\bP\to\bM$ and $\bP'\to\bM$ are semiprojective resolutions, then
for every $n\in\BZ$ there exist projective modules $Q_n$ and $Q'_n$,
such that
$$
\Ch{n}{\bP}\oplus Q'_n\cong\Ch{n}{\bP'}\oplus Q_n
$$
\end{proposition}
\begin{proof}
By \ref{semiprojective_quism}, there exists a comparison of
semiprojective resolutions $\bP'\to\bP$.  Using \ref{hom_spliting},
one can factor it as a composition of quasiisomorphisms
\[
\bP'\to\bP''\xra{\varphi}\bP
\]
where $\bP''=\bP'\oplus\Cone(\id{\bP})$ and $\varphi$ is surjective.  For
each $n\in\BZ$ we now have
\[
\Ch{n}{\bP''}=\Ch{n}{\bP'}\oplus\Ch{n}{\Cone(\id{\bP})}
\]
The complex $\Cone(\id{\bP})$ is contractible, so the module
$Q_n=\Ch{n}{\Cone(\id{\bP})}$ is projective.  On the other hand, one
has an exact sequence of complexes
\[
0\to\bQ'\to\bP''\xra{\varphi}\bP\to0
\]
Since $\varphi$ is a quasiisomorphism, we conclude that $\HH{\!}{(\bQ')}=0$. 
Because in addition $\bP$ is semiprojective, the
sequence splits yielding an isomorphism of complexes
\[
\bP''\cong\bP\oplus\bQ'
\]
It shows that $\bQ'$ is semiprojective and produces an isomorphism of modules
\[
\Ch{n}{\bP''}\cong\Ch{n}{\bP}\oplus\Ch{n}{\bQ'}
\]
By \ref{semiproj_contr}, $\bQ'$ is contractible, so the module
$Q'_n=\Ch{n}{\bQ'}$ is projective.
\end{proof}

\begin{chunk}
We define the $n^{\text{th}}$ \emph{syzygy} of a complex $\bM$ to
be $\Omega_n(\bM)=\Ch{n}{\bP}$, where $\bP\to\bM$ is a semiprojective
resolution. In view of Proposition \ref{Shanuel}, $\Omega_n\bM$ is
defined uniquely up to a projective direct summand.
\end{chunk}

\begin{proposition}[Horseshoe Lemma]
\label{semiproj_exact}
For every exact sequence of complexes
$$
0\to \bM\xra{\mu}\bM'\xra{\mu'}\bM''\to 0
$$ 
there exists a commutative diagram with exact rows
\begin{equation}
\begin{gathered}
\tag{$*$}
\label{diag*}
\xymatrixrowsep{2pc}
\xymatrixcolsep{2pc}
\xymatrix{
0\ar@{->}[r]
&\bM\ar@{->}[r]^{\mu}
&\bM'\ar@{->}[r]^{\mu'}
&\bM''\ar@{->}[r]
&0
\\
0\ar@{->}[r]
&\bP\ar@{->}[r]^{\ov\mu}\ar@{->}[u]^{\pi}
&\bP'\ar@{->}[r]^{\ov\mu'}\ar@{->}[u]^{\pi'}
&\bP''\ar@{->}[r]\ar@{->}[u]^{\pi''}
&0
}
\end{gathered}
\end{equation}
in which the columns are strict semiprojective resolutions.
\end{proposition}

\begin{proof}
Using \ref{semi_res}, choose strict semiprojective resolutions
$\pi\col\bP\to\bM$ and $\varkappa\col\bQ\to\bM'$.  Set
$\bP'=\Cone(\id{\bP})\oplus\bQ$, and note that this complex
is semiprojective.  The map $\pi'\col \bP'\to\bM'$ defined by
$\pi'(p',p,q)=\varkappa(q)$ is a composition of the canonical epimorphism
$\Cone(\id{\bP})\oplus \bQ\to\bQ$, whose kernel $\Cone(\id{\bP})$ is
contractible, with the quasiisomorphism $\varkappa$.  It follows that
$\pi'$ is a quasiisomorphism, hence $\pi'$ is a strict semiprojective
resolution.

Choose, by \ref{quism}, a morphism $\alpha\col\bP\to\bQ$ such that
$\varkappa\alpha= \mu\pi$ and define a map $\ov\mu\col\bP\to\bP'$
by $\ov\mu(p)=(0,p,\alpha(p))$. It is easy to check that this is an
injective morphism.  Set $\bP''=\Coker\ov\mu$, and form a diagram
(\ref{diag*}), where ${\ov\mu}'$ denotes the canonical surjection and
$\pi''\col\bP''\to\bM''$ is the unique morphism induced by $\pi'$
satisfying $\pi''{\ov\mu}'=\mu'\pi'$. The resulting diagram is
commutative with exact rows.

It remains to prove that $\pi''\col\bP''\to \bM''$ is a strict
semiprojective resolution.  This morphism is surjective, because
${\mu}'$ and $\pi'$ are.  It is a quasiisomorphism by the Five-Lemma,
because $\pi$ and $\pi'$ are.  For each $n$ define a map $\beta_n\col
P'_n\lra P_{n-1}\oplus Q_n$ by $\beta(p',p,q)=(p',q-\alpha(p))$.  These
maps are clearly $R$-linear, surjective, with $\Ker\beta_n=\im\ov\mu_n$,
so they produce isomorphisms $P''_n\cong P_{n-1}\oplus Q_n$.  Thus,
every $R$-module $P''_n$ is projective.  We now know that the lower
row of the diagram splits as a sequence of graded $R$-modules, so each
complex $\bN$ induces an exact sequence
\[
\xymatrixrowsep{2pc}
\xymatrixcolsep{2pc}
\xymatrix{
0\ar@{->}[r]
&{\Hom R{\bP''}\bN}\ar@{->}[r]
&{\Hom R{\bP'}\bN}\ar@{->}[r]
&{\Hom R{\bP}\bN}\ar@{->}[r]
&0
}
\]
As $\bP'$ and $\bP$ are semiprojective, $\HH{\!}{(\bN)}=0$ implies
$\HH{\!}{(\Hom R{\bP'}\bN)}=0$ and $\HH{\!}{(\Hom R{\bP}\bN)}=0$
by \ref{semiprojective}.  From the homology exact sequence one now
gets $\HH{\!}{(\Hom R{\bP''}\bN)}=0$, so $\bP''$ is semiprojective  by
\ref{semiprojective}.  Summing up the preceding discussion, we see
that $\pi''\col\bP''\to\bM''$ is a strict semiprojective resolution.
\end{proof}

\subsection{Projective dimension} 

In \cite{af} the \emph{projective dimension} of $\bM$ is defined to be
the infimum of those $p\in\BZ$ for which there exists a semiprojective resolution
$\bP\to\bM$ with $P_i=0$ for all $i>p$; when no such $p$ exits, $\bM$
is said to have infinite projective dimension.  The projective dimension
of $\bM$ is denoted $\pd R\bM$.  When $\bM$ is a module, this number is
equal to the classical invariant.

\begin{chunk}\cite[(2.4.P)]{af}
\label{projdim}
For each integer $p$ the following conditions are equivalent. 
\begin{enumerate}[\rm\quad(i)]
\item
$\pd R\bM\le p$.
\item
$\sup\HH{\!}{(\bM)}\leq p$ and there exists a semiprojective resolution
$\bP\to\bM$ such that the module $\Ch{p}{\bP}$ is projective.
\item[(ii$'$)]
$\sup\HH{\!}{(\bM)}\leq p$ and for every semiprojective resolution
$\bP'\to\bM$ the module $\Ch{p}{\bP'}$ is projective. 
\item
$\sup\HH{\!}{(\bM)}\leq p$ and $\Omega_{p}{(\bM)}$ is projective.
\end{enumerate}
\end{chunk}

\begin{chunk}
\cite[(2.4.P)]{af}
\label{projdim_ext} 
There are equalities
\begin{align*}
\pd{R}{\bM}&=\sup\left\{n\in\BZ\,\left|\,\begin{gathered}
\Ext{R}{n}{\bM}{N}\not=0 \\
\text{for some module}\  N\end{gathered}\right.\right\}\\
&=\sup\left\{n\in\BZ\,\left|\,\begin{gathered}
\Ext R{n-\inf\HH{\!}{(\bN)}}\bM\bN\not=0\\
\text{for some complex}\ \bN\  \text{with}\ \inf\HH{\!}{(\bN)}>-\infty
\end{gathered}\right.\right\}.
\end{align*}
Thus,  projective dimension of complexes can be computed in terms of
vanishing of appropriate Ext groups, as is the case for modules.
\end{chunk}

An immediate  consequence of Proposition \ref{semiproj_exact} and
 \ref{projdim_ext} is the following. 

\begin{chunk}
\label{pd_ineq}
Let  $0\to \bM\to \bM'\to \bM''\to 0$ be an exact sequence of
complexes.
\begin{enumerate}[\rm\quad(1)]
\item If two complexes have finite projective dimension,
then so does the third.
\item There is an inequality
$$
\pd R{\bM'}\le\max(\pd R{\bM},\pd R{\bM''})
$$
and equality holds, except possibly when $\pd R{\bM''}=\pd
R{\bM}+1$.
\end{enumerate} 
\end{chunk}

\section{Totally acyclic complexes}
\label{Totally acyclic}

In this section we first review the notions of totally acyclic complex
and complete resolution, which are at the basis of the definition of
Gorenstein projective dimension.  Second, we collect properties of modules
of finite Gorenstein projective dimension of modules that will be used
to define and study Gorenstein projective dimension and Tate cohomology
for complexes.  Finally, we compare the notion of Gorenstein projective
dimension with the earlier notion of Gorenstein dimension, which is 
defined only for finite modules over noetherian rings.

\subsection{Total acyclicity}
The following notion is fundamental to our investigation.

\begin{chunk}
\label{Totally_acyclic}
A complex of $R$-modules $\bT$ is said to be
\emph{totally acyclic} if the following conditions are satisfied.
\begin{enumerate}[\rm\quad (1)]
\item $T_n$
is projective for every $n\in\BZ$.
\item $\bT$ is  exact.
\item $\Hom{R}{\bT}{Q}$
is exact for every projective $R$-module $Q$.
\end{enumerate} 
\end{chunk}

This meaning of ``totally acyclic'' is wider than the one in \cite{am},
where it refers to an exact complex $\bT$ of \emph{finite} projective
modules such that $\Hom{R}{\bT}{R}$ is exact; when $R$ is noetherian
the new notion subsumes the earlier one, by \cite[(2.4)]{am}.

We shall use the following result of Cornick and Kropholler.

\begin{chunk}
\cite[(2.4)]{ck}
\label{extension}
Let $\bT$ be a totally acyclic complex. If $\bQ$ is a complex of
projectives and $n$ is an integer, then any chain map
$\ul{\phi}\col\bT_{\geq n}\to\bQ_{\geq n}$ can be extended to a chain map
$\phi\col\bT\to\bQ$ such that $\phi_{\geq n}=\ul{\phi}$. Every chain map $\phi$  with this property is defined unique up to homotopy.
\end{chunk}

As a consequence we obtain the following.
\begin{chunk}
\label{total_vanishing}
Let $\bT$ be a totally acyclic complex. If  $\bQ$ is a bounded above
complex of projectives, then
$$
\HH{\!}{(\Hom{R}{\bT}{\bQ})}=0.
$$
\end{chunk}
Indeed, set $\sup\bQ=s$. If $\alpha\col\bT\to\bQ$ is a chain map of degree $n$, then $\alpha_{\geq s-n+1}$ is the zero map. The preceding result shows that $\alpha$ is homotopic to $0$, that is $\HH{n}{(\Hom{R}{\bT}{\bQ})}=0$.

\subsection{Complete resolutions}
\label{Complete resolutions} 
A notion of complete resolution was initially introduced, for finite modules over finite groups, by Tate \cite[Ch. XII]{ce}. Cornick and Kropholler used it in \cite{ck}, in a modified form; more precisely, they added the condition (3) in the definition of a  totally  acyclic complex \ref{Totally_acyclic}. Avramov and Martsinkovsky \cite{am} incorporated a comparison morphism to a projective resolution, in the notion of complete resolution for finite modules over noetherian rings. Using this last perspective, we extend the notion of complete resolutions to complexes. 

Let $\bM$ be a complex of $R$-modules.

\begin{chunk}
\label{complete resolution}
A \emph{complete resolution} of $\bM$ is a diagram of morphisms of
complexes 
$$
\bT\xrightarrow{\tau}\bP\xrightarrow{\pi}\bM
$$
where
$\pi\col\bP\to\bM$ is a semiprojective resolution, $\bT$ is a totally
acyclic complex and $\tau_i$ is bijective for all $i\gg 0$. A complete resolution is said to be  \emph{surjective} if $\tau_i$ is
surjective for all $i\in\BZ$.
\end{chunk}

The next result shows that from any complete resolution one can get
a surjective complete resolution.  The proof given in \cite[(3.7)]{am}
extends easily and is omitted.

\begin{chunk}
\label{surjective_complete}
Let  $\bT\stackrel\tau\to\bP\xrightarrow{\pi}\bM$ be a complete
resolution. If $g$ is an integer such that  $\tau_{i}$ is
bijective for all $i\ge g$, then there exists a complete resolution
$\bT'\stackrel{\tau'}\to\bP\xrightarrow{\pi}\bM$ where $\tau'$ is a
surjective morphism, such that $\tau'_{i}$ bijective for all $i\ge g$, and
a homotopy equivalence $\alpha\col\bT\to\bT'$ such that $\tau=\tau'\alpha$
and $\alpha_i=\id{T_i}$ for all $i\geq g$.  
\end{chunk}

\begin{chunk}
\label{ker_complete}
Let $\bT\xra{\tau}\bP\xra{\pi}\bM$ be a surjective complete
resolution and set $\bL=\Ker\tau$.  The short exact
sequence of complexes 
$$
0\to\bL\xra{\varkappa}\bT\xra{\tau}\bP\to0
$$
splits as a sequence of graded modules since $\bP$ is a complex
of projectives.  Let  $\lambda\col\bP\to\Sigma\bL$ be the morphism
$\theta^{\bP}$ constructed in \ref{theta_maps}. By \ref{theta_morphism},
we have $\HH{\!}{(\lambda)}=\eth$, where $\eth$ is the connecting
homomorphism of the homology exact sequence associated to
the short exact sequence above. Since $\bT$ is exact, $\lambda$ is a
quasiisomorphism.
\end{chunk}

When $\bM$ and $\bM'$ are finite modules over noetherian rings, the  following result is contained in the statement and the proof  of \cite[(5.3)]{am}. The arguments carry over verbatim, so they are omitted.
\begin{proposition}
\label{Tate_lifting}
Let $\bT\xra\tau\bP\xra\pi\bM$ and $\bT'\xra{\tau'}\bP'
\xra{\pi'}\bM'$ be complete resolutions. For each morphism of
complexes $\mu\col \bM\to\bM'$ there exists a unique up to homotopy
morphism $\ov\mu$, making the right hand square of the diagram
\[
\xymatrixrowsep{2pc}
\xymatrixcolsep{2pc}
\xymatrix{
&\bT\ar@{->}[d]^{\wh\mu}\ar@{->}[r]^{\tau}
&\bP\ar@{->}[d]^{\ov\mu}\ar@{->}[r]^{\pi}
&\bM\ar@{->}[d]^{\mu}
\\
&\bT'\ar@{->}[r]^{\tau'}
&\bP'\ar@{->}[r]^{\pi'}
&\bM'
}
\]
commute up to homotopy, and for each choice of $\ov\mu$ there exists a unique up
to homotopy morphism $\wh\mu$, making the left hand square commute
up to homotopy. Moreover, if $\pi'$ and $\tau'$ are surjective, then $\ov\mu$ and respectively $\wh\mu$ can be chosen such that the right hand square and respectively the left hand side square commute.

If $\mu=\id{\bM}$, then $\ov\mu$ and $\wh\mu$ are homotopy
equivalences.\QED
\end{proposition}

\subsection{Gorenstein projective modules}
\label{The case of modules}
A notion of Gorenstein projective dimension, for left modules over
associative rings, is due to Enochs and Jenda \cite{ej1}. It extends
the notion of Gorenstein dimension, introduced by Auslander and
Bridger \cite{ab}, for finite modules over commutative noetherian
rings. Enochs and Jenda study it for left or right modules over coherent or
n-Gorenstein rings and Holm \cite{hh} extends their results to arbitrary modules over associative rings. 
In this subsection we recall some  definitions and give some basic results that are needed in this paper.

\begin{chunk}
\label{G-projective}
An $R$-module $G$ is \emph{Gorenstein projective} if there exists a totally acyclic
complex $\bT$, as defined
in \ref{Totally_acyclic}, with $\Ch{0}{\bT}=G$.
\end{chunk}

\begin{chunk}\cite[(2.20)]{hh}
\label{Gp_Ext}
If $G$ is Gorenstein projective and $Q$ is projective, then
$$
\Ext{R}{i}{G}{Q}=0\qquad \text{for all}\qquad i\geq 1.
$$ 
Indeed, by definition there exists a totally acyclic complex $\bT$
with $\Ch{0}{\bT}=G$.  Therefore, $\bT_{\ge 0}\to G$ is a projective
resolution, so for all $i>0$ one has
$$
\Ext{R}{i}{G}{Q}=\Hh{i}{(\Hom{R}{\bT_{\ge 0}}{Q})}=\Hh{i}{(\Hom{R}{\bT}{Q})}=0.
$$
\end{chunk}

\begin{lemma}
\label{syz_exact} 
For an exact complex of projectives $\bT$ the following are
equivalent.
\begin{enumerate}[\rm\quad(i)]
\item
$\bT$ is totally acyclic.
\item
$\Ch{i}{\bT}$ is Gorenstein projective for all $i\in\BZ$.
\item
$\Ch{i}{\bT}$ is Gorenstein projective for infinitely many $i\le 0$.
\end{enumerate}
\end{lemma}

\noindent\emph{Proof.} 
The implications (i)$\implies$(ii)$\implies$(iii) are clear, so we only
argue (iii)$\implies$(i).  Let $Q$ be a projective $R$-module and fix
$n\in\BZ$. We need to show $\Hh{n}{(\Hom{R}{\bT}{Q})}=0$. Choose an
integer $m\geq 1$ so that  $n-m$ is small enough for $\Ch{n-m}{\bT}$
to be Gorenstein projective. The conclusion follows via \ref{Gp_Ext}:
\begin{xxalignat}{3}
&\hphantom{\square}
&\Hh{n}{(\Hom{R}{\bT}{Q})}&=\Ext{R}{m}{\Ch{n-m}{\bT}}{Q}=0.
&&\square
\end{xxalignat}

\begin{lemma}
\label{res_tot}
If $G$ is Gorenstein projective and $\bP\to G$ is a projective
resolution, then there exists a totally acyclic complex $\bT$ such
that $\bT_{\geq 0}=\bP$. 
\end{lemma}

\begin{proof} 
By definition, there is a totally acyclic complex $\bT'$ with
$\Ch{0}{\bT'}=G$. Set
 \begin{alignat*}{3}
T_i = \begin{cases}T'_i &\text{for}\ i< 0\,;\\ 
                   P_i&\text{for}\  i\geq 0\,;
      \end{cases}
\qquad
&\text{and}&
\qquad
\partial_i^{\bT}=\begin{cases}\partial_i^{\bT'}&\text{for}\ i<0\,;\\
                              \beta\alpha&\text{for}\ i=0\,;\\
                              \partial_i^{\bP}&\text{for}\  i>0\,,
               \end{cases}
\end{alignat*}
where $\alpha\col P_0\to G$ and $\beta\col G\to T'_{-1}$  are the
canonical maps. The complex $\bT$ is exact, $\Ch{0}{\bT}=G$ and
$\bT_{<0}=\bT'_{<0}$, so $\bT$ is totally
acyclic by Lemma \ref{syz_exact}.
 \end{proof}

A basic property of Gorenstein projective modules is proved by Holm
\cite[(2.5)]{hh}.

\begin{chunk}
\label{resolving}
The class of Gorenstein projective modules has the following properties.
\begin{enumerate}[\rm\quad(1)]
\item
In every exact sequence $0\to G\to G'\to G''\to 0$ of $R$-modules with
$G''$ Gorenstein projective, the module $G$ is Gorenstein projective if and
only if the module $G'$ is Gorenstein projective.
\item
When $G'\cong G\oplus G''$ the module $G'$ is Gorenstein projective
if and only if both modules $G$ and $G''$ are Gorenstein projectives.
\item 
Every projective module is Gorenstein projective.
\end{enumerate}
\end{chunk}

Let $M$ be a left  $R$-module.

\begin{chunk}
\label{Gpd for modules}
Following Enochs and Jenda \cite{ej1},
we say that $M$ has \emph{Gorenstein projective dimension} $g$, and write $\Gpd R M=g$, if there exists
an exact sequence 
$$
0\rightarrow G_g\rightarrow G_{g-1}\rightarrow\cdots\rightarrow
G_1\rightarrow G_{0}\rightarrow M\rightarrow 0
$$
where $G_i$ is Gorenstein projective for $i=0,\dots,g$, and there 
is no shorter exact sequence with this property.
\end{chunk}

The next three results are due to Holm.

\begin{chunk}\cite[(2.7)]{hh}
\label{ext1}
One has $\Gpd R M\le n$ if in one, and only if in
every exact sequence 
$$
0\rightarrow K_n\rightarrow G'_{n-1}\rightarrow
\cdots\rightarrow G'_1\rightarrow G'_{0}\rightarrow M\rightarrow 0
$$
where all $G'_i$ are Gorenstein projective, the module $K_n$ 
is Gorenstein projective.
\end{chunk}

\begin{chunk}\cite[(2.27)]{hh}
\label{ext2}
One has $\Gpd R M\le\pd R M$, with equality when $\pd R M$ is finite.
\end{chunk}

\begin{chunk}\cite[(2.19)]{hh}
\label{direct}
For every family of $R$-modules $(M_i)_{i\in I}$ one has
$$
\Gpd R {\big(\coprod_{i\in I}M_i\big)}=\sup_{i\in I}\{\Gpd R{M_i}\}
$$
\end{chunk}


\subsection{Gorenstein dimension}
\label{Gorenstein dimension}
In this subsection we assume that the ring $R$ is left and right
noetherian and that $M$ is a finite $R$-module.

Auslander and Bridger \cite{ab} define the 
\emph{Gorenstein dimension} $\Gdim RM$, as
the supremum of those integers $n$ for which there exists an exact
sequence
$$
0\rightarrow G_n\rightarrow\cdots\rightarrow G_1\rightarrow
G_{0}\rightarrow M\rightarrow 0
$$
where each $G_i$ is a finite reflexive $R$-module
satisfying 
$$
\Ext {R}{n}{G_i}R=0=\Ext {R}{n}{G_i^*}R \qquad\text{for all}\qquad n\ge1\,.
$$

\begin{chunk}
\label{agree}
The two notions of Gorenstein dimension agree:
$$
\Gdim RM=\Gpd RM
$$
Indeed, $\Gdim RM=0$ if and only $\Gpd RM=0$ by \cite[(4.2.6)]{cl};
one deduces from here that equality holds in general using \ref{ext1}. 
\end{chunk}

Let $\idim RR$ (respectively, $\idim{R^{\rm o}}R$) denote the injective
dimension of $R$ as a left (respectively, right) $R$-module.
 
\begin{chunk}
\label{Gorenstein_Gpd}  
Let $n$ be a non-negative integer.  The following conditions are equivalent.
\begin{enumerate}[\rm\quad(i)]
\item
$\idim RR\le n$ and $\idim{R^{\rm o}}R\le n$.
\item
Every $R$-module $M$ satisfies $\Gpd RM\le n$.
\item
Every finite $R$-module $M$ satisfies $\Gdim RM\le n$.
\end{enumerate}

Indeed, (i) implies (ii) by \cite[(10.2.14)]{ej2}, (ii) implies
(iii) by \ref{agree}, and (iii) implies (i) by \cite[(3.2)]{am}.
 \end{chunk}

Gorenstein dimension has been extensively studied over commutative
rings.  We quote some results that are used in this paper.

\begin{chunk}
\label{ABformula}
Let $R$ be a commutative ring and $M$ a finite $R$-module.
\begin{enumerate}[\rm\quad(1)]
\item\cite[(4.15)]{ab}
If $\bsp$ is a prime ideal of $R$, then
$$
\Gdim{R_\bsp}{M_\bsp}\le\Gdim RM
$$
\item\cite[(4.13.b)]{ab}
If $R$ is local and $\Gdim RM$ is finite, then
$$\Gdim RM+\depth R M=\depth{}R$$
\end{enumerate}
\end{chunk}

We use the term \emph{Gorenstein ring} to denote a noetherian commutative
ring $R$, such that the local ring $R_\bsm$ has finite injective dimension
as a module over itself for every maximal ideal $\bsm$ of $R$.

\begin{chunk}
\label{Gring}
If $R$ is Gorenstein,
then $\idim RM=\dim R$, see \cite[(3.1.17)]{bh} and \cite[(3.2.20)]{ej2}.
\end{chunk}

The next result, due to Auslander and Bridger \cite[(4.20)]{ab} and Goto
\cite[Corollary 2]{go}, gives the reason for the name \emph{Gorenstein
dimension}.

\begin{chunk}
\label{Goto}
If $R$ is a commutative noetherian ring, then $R$ is Gorenstein if and
only if every finite module has finite Gorenstein dimension.
\end{chunk}

\section{Gorenstein projective dimension}
\label{Gorenstein projective}

In this section, first we define a notion of Gorenstein projective dimension
for complexes which extends naturally the notions of Gorenstein
(projective) dimensions, previously defined by several
authors, see Corollary \ref{old_def} and Remark \ref{old_new}. Second, we prove our main results, Theorems \ref{complete},
\ref{complete_Gpd} and respectively \ref{Gpd_ineq} which are parallel
to the results on projective dimension \ref{projdim},
\ref{projdim_ext} and respectively \ref{pd_ineq}, and Theorem
\ref{Gpd_Gorenstein} which generalizes \ref{Goto}.

Let $\bM$ be a complex of left $R$-modules. 

\begin{5definition}
\label{Gpd for complexes}
The \emph{Gorenstein projective dimension} of $\bM$ is defined by

$$\Gpd R\bM=\inf\left\{g\in\BZ\left |\begin{gathered} \bT\xra\tau\bP\xra\pi\bM\\ \text{is a complete resolution with}\\  \tau_{i}\col T_{i}\to
P_{i}\  \text{bijective for each}\ i\ge g \end{gathered}\right.\right\}$$
 \end{5definition}

The next two assertions follow from the definition.
 
\begin{5chunk}
\label{Gpd_exact} 
The complex  $\bM$ is  exact if and only if $\Gpd R \bM=-\infty$. 
\end{5chunk}

\begin{5chunk}
\label{shiftt}
For every $i\in\BZ$,  one has $\Gpd R \Sigma^i\bM=\Gpd R\bM +i$.  
\end{5chunk}

The next theorem can be viewed as an equivalent definition of $\Gpd
R\bM$, expressed in terms of an arbitrary semiprojective resolution
of the complex $\bM$.
   
\begin{5theorem}
\label{complete}
For each integer $g$ the following conditions are equivalent.
\begin{enumerate}[\rm\quad(i)]
\item
$\Gpd R\bM\le g$.
\item   
$\sup \HH{\!}{(\bM)}\leq g$ and there exists a semiprojective resolution
$\bP\to\bM$ such that the module $\Ch{g}{\bP}$ is Gorenstein projective.
\item[\rm (ii$'$)]
$\sup \HH{\!}{(\bM)}\leq g$ and for every semiprojective resolution
$\bP'\to\bM$ the module $\Ch{g}{\bP'}$ is Gorenstein projective.
\item 
$\sup \HH{\!}{(\bM)}\leq g$ and $\Omega_g(\bM)$ is Gorenstein projective.
\item 
For every semiprojective resolution $\bP'\to\bM$ there exists a surjective
complete resolution $\bT'\xra{\tau'}\bP'\to\bM$ such that
$\tau_i'=\id{T'_i}$ for all $i\geq g$.
\end{enumerate}
\end{5theorem}

\begin{proof}
(i)$\implies$(ii).
By hypothesis, there is a complete resolution
$\bT\xra{\tau}\bP\xra\pi\bM$ such that $\tau_{\geq g}\col \bT_{\ge
g}\to \bP_{\ge g}$ is an isomorphism of complexes.  This implies
isomorphisms  $\HH{i}{(\bM)}\cong\HH{i}{(\bP)}$ for all $i\in\BZ$, $\HH{i}{(\bP)}\cong\HH{i}{(\bT)}$ for all $i>g$, and
$\Ch{g}{\bP}\cong\Ch{g}{\bT}$. Since the complex $\bT$ is totally
acyclic, $\HH{i}{(\bT)}$ vanishes for each $i\in\BZ$ and $\Ch{g}{\bT}$
is Gorenstein projective. 
 
(ii)$\implies$(ii$'$). 
Let $\bP'\xrightarrow{\pi'}\bM$ be a semiprojective resolution.  By Proposition \ref{Shanuel}, there exist projective modules $Q_g$
and $Q_g'$ such that $\Ch{g}{\bP}\oplus Q_g'\cong\Ch{g}{\bP'}\oplus
Q_g$. Thus, $\Ch{g}{\bP'}$ is Gorenstein projective by \ref{resolving}.

(ii$'$)$\implies$(iii) is clear.

(iii)$\implies$(iv).
Let $\bP'\to\bM$ be a semiprojective resolution with
$\Omega_g{(\bM)}=\Ch{g}{\bP'}$ Gorenstein projective and
$\HH{i}{(\bP')}=0$ for all $i>g$. Then $\Sigma^{-g}\bP'_{\geq
g}\to\Ch{g}{\bP'}$ is a projective resolution. Lemma \ref{res_tot},
yields a totally acyclic complex $\bT''$ with $\bT''_{\geq
g}=\bP'_{\geq g}$.  Applying \ref{extension}, we obtain a complete
resolution $\bT''\xrightarrow{\tau''}\bP'\xrightarrow{}\bM$ with
$\tau''_i=\id{T''_i}$ for all $i\geq g$, and $\Ch{g}{\bT''}\cong
\Ch{g}{\bP'}$. From \ref{surjective_complete} we get a surjective complete
resolution $\bT'\xra{\tau'}\bP'\xra{\pi'}\bM$ with the desired properties.

(iv)$\implies$(i) is clear.\end{proof}

\begin{5corollary}
\label{coproducts}
For every family of complexes of $R$-modules $(\bM_i)_{i\in I}$ one has
$$
\Gpd R {\big(\coprod_{i\in I}\bM_i\big)}=\sup_{i\in I}\{\Gpd R{\bM_i}\}
$$
\end{5corollary}

\begin{proof}
Choose for each $i\in I$ a semiprojective resolution $\bP_i\to\bM_i$
and set $\bP=\coprod_{i\in I}\bP_i$.  As $\bP\to
\coprod_{i\in I}\bM_i$ is a semiprojective resolution and
$\Ch{n}{\bP}=\coprod_{i\in I}\Ch{n}{\bP}$ for each $n\in\BZ$,
the assertion follows from the theorem and \ref{direct}.
\end{proof}

We reconcile the two notions of Gorenstein projective dimension for
modules.

\begin{5corollary}
\label{old_def}
If $M$ is an $R$-module and  $g$ (respectively, $g'$) denote the Gorenstein projective dimension
of $M$ in the sense of \ref{Gpd for modules} (respectively, of \ref{Gpd
for complexes}), then $g=g'$.
 \end{5corollary}

\begin{proof}
Choose a classical projective resolution $\bP\to
M$, and let $h$ be the smallest integer $n$ such that $\Ch{n}{\bP}$ is
Gorenstein projective.  As each $P_i$ is Gorenstein projective, \ref{ext1}
implies $h=g$.  On the other hand, as $\bP\to M$ is a semiprojective
resolution  by \ref{proj}, Theorem \ref{complete} yields $h=g'$.
 \end{proof}

Another consequence of the Theorem \ref{complete} is an extension to complexes of \ref{ext2}.

\begin{5theorem}
\label{Gpd_pd}
One has $\Gpd R\bM\le\pd R\bM$, with equality if $\pd R\bM$ is finite.
\end{5theorem}

\begin{proof}
Set $p=\pd R\bM$ and $g=\Gpd R\bM$. There is nothing to
prove if $p=\infty$.  If $p=-\infty$, then, by definition
$\HH{\!}{(\bM)}=0$. Therefore, by \ref{Gpd_exact},
$\Gpd{R}{\bM}=-\infty$. We may assume that $p$ is finite.

By \ref{projdim},  $\HH{i}{(\bM)}=0$ for $i>p$ and  there exists a
semiprojective resolution $\bP\to\bM$ in which the module $\Ch{p}{\bP}$ is
projective, hence Gorenstein projective. Applying Theorem \ref{complete},
we obtain $g\leq p$.

Suppose $g<p$.  By definition, there exists a semiprojective resolution
$\bP'\to\bM$ such that $P'_i=0$ for all $i>p$. By Theorem \ref{complete},
$\Ch{g}{\bP'}$ is Gorenstein projective and $\HH{i}{(\bM)}=0$ for all
$i>g$,  in particular $\HH{i}{(\bP')}=0$ for all $i>g$. Therefore,
there exists an exact sequence of $R$-modules
\[
0\rightarrow P'_p\rightarrow\cdots\rightarrow P'_g\rightarrow \Ch{g}{\bP'}
\rightarrow 0
\]
Thus, the Gorenstein projective module $\Ch{g}{\bP'}$ has finite
projective dimension, so it is projective by \ref{ext2}. Now
\ref{projdim} implies $p\le g$, which contradicts our assumption.
Hence $p=g$, as desired.
\end{proof}

The next theorem extends to complexes a result  of Holm,
\cite[(2.20)]{hh}.  While formally similar to the characterization of
projective dimension in terms of vanishing of appropriate Ext groups, see
\ref{projdim_ext}, it differs in a significant aspect: it is restricted to complexes of \emph{finite} Gorenstein projective dimension.

\begin{5theorem}
\label{complete_Gpd}
If $\Gpd{R}{\bM}=g<\infty$, then there are equalities
\begin{align*}
\Gpd{R}{\bM}
&=\sup\left \{n\in\BZ\,\left|\,\begin{gathered}
\Ext{R}{n}{\bM}{Q}\not=0\\
\text{for some projective module}\ Q
\end{gathered}\right \}\right.\\
&=\sup\left \{n\in\BZ\,\left|\,\begin{gathered}
\Ext{R}{n-\inf\HH{\!}{(\bN)}}{\bM}{\bN}\ne0\\
\text{for some complex}\ \bN \text{with}\\
 \pd{R}{\bN}<\infty\ \text{and}\ \inf\HH{\!}{(\bN)>-\infty}
\end{gathered}\right \}\right..
\end{align*}
\end{5theorem}

\begin{proof}
When $\Gpd R\bM=-\infty$ we use \ref{Gpd_exact}. Suppose
that $g$ is finite.  By definition, there is a complete resolution
$\bT\xrightarrow{\tau}\bP\xrightarrow{\pi}\bM$ such that $\tau_{\geq
g}\col\bT_{\geq g}\to\bP_{\geq g}$ is an isomorphism.  Let $h$
(respectively, $h'$) denote the number on the right hand side of the first
(respectively, second) equality of the theorem.

First, we show  $g\ge h'$.  By \ref{semi_res} and \ref{projdim}
choose a  semiprojective resolution  $\bQ\to\bN$ with
$\inf\bQ=\inf\HH{\!}{(\bN)}$ and $\sup\bQ<\infty$.  For all
$i>g-\inf\HH{\!}{(\bN)}$ we have
\begin{align*} 
\Ext{R}{i}{\bM}{\bN}&=\Hh{i}{(\Hom{R}{\bP}{\bN})}\\
                    &\cong\Hh{i}{(\Hom{R}{\bP}{\bQ})}\\
                    &=\Hh{i}{(\Hom{R}{\bT}{\bQ})}\\
                    &=0.
\end{align*}
Indeed, the first equality holds by definition; the isomorphism is
due to the semiprojectivity of $\bP$; for the third equality we use
the isomorphism $\bP_{\geq g}\cong T_{\geq g}$ and the choice of $i$;
the last equality follows from \ref{total_vanishing}.

The inequality $h'\ge h$ is obvious, so it remains to prove $h\ge g$.

By way of contradiction, assume $h<g$.  The sequence
\begin{alignat*}{4}
\Hom{R}{P_{g-1}}{T_{g-1}}&\xra{\Hom{R}{\partial^\bP_{g}}{T_{g-1}}}&&\Hom{R}{P_{g}}{T_{g-1}}\\
                         &\xra{\Hom{R}{\partial^\bP_{g+1}}{T_{g-1}}}&&\Hom{R}{P_{g+1}}{T_{g-1}}
\end{alignat*}
is then exact.  Let $\alpha\col P_g\cong T_g \to T_{g-1}$ be the canonical map. Since $\Hom{R}{\partial^\bP_{g+1}}{T_{g-1}}(\alpha)=0$, there exists a map $\beta\col P_{g-1}\to T_{g-1}$ such that  $\beta\partial_{g}^\bP=\alpha$. It is easy to see that the pair  $(\bT',\partial^{\bT'})$ given by 
\begin{align*}
T'_i=\begin{cases}
            T_i&\text{for}\ i\not={g-1}\,;\\
            P_{g-1}&\text{for}\ i=g-1\,;
     \end{cases}
\qquad
& \text{and}
\qquad
\partial^{\bT'}_i=\begin{cases} 
                        \partial^\bT_i&\text{for}\ i\not=g,\ g-1\,;\\
                        \tau_{g-1}\partial^\bT_g&\text{for}\ i=g\,;\\
                        \partial^\bT_{g-1}\beta&\text{for}\ i={g-1}\,,
                  \end{cases}
\end{align*}
is an exact complex. Since each module $T_i'$ is projective and $\Ch{i}{\bT'}$ is Gorenstein projective for all  $i\ll 0$, $\bT'$ is totally acyclic  by \ref{syz_exact}. Thus, we obtain a complete resolution $\bT'\xra{\tau'}\bP\xra{\pi}\bM$ with $\tau'_i$ bijective for each $i\geq g-1$. Therefore, we get $\Gpd R\bM\leq g-1$; this is a contradiction, hence $h=g$. 
\end{proof}

The next result extends \cite[(4.9.2)]{am} to complexes 
and gives a parallel to \ref{pd_ineq}.

\begin{5theorem}
\label{Gpd_ineq}
Let  $0\to \bM\to \bM'\to \bM''\to 0$ be an exact sequence of
complexes.
\begin{enumerate}[\rm\quad(1)]
\item If two complexes have finite Gorenstein projective dimension,
then so does the third.
\item There is an inequality
$$
\Gpd R{\bM'}\le\max(\Gpd R{\bM},\Gpd R{\bM''})
$$
and equality holds, except possibly when $\Gpd R{\bM''}=\Gpd
R{\bM}+1$.
\end{enumerate} 
\end{5theorem}

\begin{proof}
(1). By  Proposition \ref{semiproj_exact} there exists a commutative
diagram
 \[
\xymatrixrowsep{2pc}
\xymatrixcolsep{2pc}
\xymatrix{
0\ar@{->}[r]
&\bM\ar@{->}[r]^{}
&\bM'\ar@{->}[r]^{}
&\bM''\ar@{->}[r]
&0
\\
0\ar@{->}[r]
&\bP\ar@{->}[r]^{}\ar@{->}[u]^{}
&\bP'\ar@{->}[r]^{}\ar@{->}[u]^{}
&\bP''\ar@{->}[r]\ar@{->}[u]^{}
&0
}
\]
where the columns are semiprojective resolutions and the bottom short
sequence is exact.  The exact sequences
of homology groups shows that if two of the
semiprojective complexes have bounded above cohomology, then so has
the third. Using Theorem \ref{complete}, choose $n$ so that
$$
\HH{i}{(\bP)}=\HH{i}{(\bP')}=\HH{i}{(\bP'')}=0\quad\text{for all}\quad i\geq n.
$$
For all $i\geq n$ the sequence of modules
$$
0\to\Ch{i}{\bP}\to\Ch{i}{\bP'}\to\Ch{i}{\bP''}\to0
$$
is then exact.  If $\Ch{i}{\bP}$ (respectively $\Ch{i}{\bP'}$) and
$\Ch{i}{\bP''}$ are Gorenstein projective, then \ref{resolving}
shows that so is $\Ch{i}{\bP'}$ (respectively $\Ch{i}{\bP}$).
If $\Ch{i}{\bP}$ and $\Ch{i}{\bP'}$ are Gorenstein projectives, then
$\Gpd R{\Ch{i}{\bP''}}\le1$; hence  by \ref{ext1}, $\Ch{i+1}{\bP''}$
is Gorenstein projective. Applying Theorem \ref{complete}, we obtain
the desired conclusion.

(2). By (1) we may assume that all complexes have finite Gorenstein
projective dimension.  The conclusion then follows by considering the
exact sequences of groups $\Ext{R}{}{-}{Q}$, induced by the short exact
sequence from hypothesis, where  $Q$ varies over all projective modules,
and applying Theorem \ref{complete_Gpd}.  \end{proof}

The next result extends  \ref{Gorenstein_Gpd}.  It characterizes Gorenstein
rings in terms of the finiteness of the Gorenstein projective dimension of
complexes.

\begin{5theorem}
\label{Gpd_Gorenstein}
If $R$ is a left and right noetherian ring and  $n$ is a
 non-negative integer, then the following are equivalent.
\begin{enumerate}[\rm\quad(i)]
\item
$\idim RR\le n$ and $\idim{R^{\rm o}}R\le n$.
\item
Every complex $\bM$ of left $R$-modules or right $R$-modules satisfies
$$ \Gpd R\bM\leq\sup\HH{\!}{(\bM)}+n.$$ 
\end{enumerate}
\end{5theorem}

\begin{proof} 
(ii)$\implies$(i) holds by \ref{Gorenstein_Gpd}.

To prove (i)$\implies$(ii) it is enough to consider left modules
since all the results used in the argument hold also for right modules. 

We may assume $\sup\HH{\!}{(\bM)}=t<\infty$.  Set $g=t+n$ and choose a
semiprojective resolution $\bP\to\bM$. By \ref{Gorenstein_Gpd}, every
module has Gorenstein projective dimension at most $n$, in particular
$\Gpd R{\Ch{t}{\bP}}\leq n$. From the exact sequence
$$
0\to\Ch{g}{\bP}\to P_{g-1}\to\cdots\to P_t\to \Ch{t}{\bP}\to 0
$$ 
one concludes that the module $\Ch{g}{\bP}$ is Gorenstein projective,
see \ref{ext1}. Therefore, Theorem \ref{complete} yields $\Gpd R\bM\leq
g$.
 \end{proof}

Another variation on the same theme was pointed out by Avramov
and Iyengar.
 
\begin{5theorem}
\label{Gpd_Gorenstein_com}
Let $R$ be a commutative noetherian ring.
\begin{enumerate}[\rm\quad(1)]
\item
$R$ is Gorenstein and $\dim R$ is finite if and only if $\Gpd
R\bM$ is finite for every complex of $R$-modules with
$\HH{\!}{(\bM)}$ bounded above.
\item
$R$ is Gorenstein if and only if $\Gpd R\bM$ is finite for every
complex $\bM$ of $R$-modules with $\HH{\!}{(\bM)}$ bounded and
$\HH i{(\bM)}$ finitely generated for each $i\in\BZ$.
\end{enumerate}
\end{5theorem}

\begin{proof}
(1).  If $R$ is Gorenstein and $\dim R$ is finite, then $\idim RR$ is
finite by \ref{Gring}, so $\Gpd R\bM<\infty$ for every complex with
$\sup\HH{\!}{(\bM)}<\infty$ by Theorem \ref{Gpd_Gorenstein}.

Assume now that $\Gpd RM$ is finite for every $R$-module $M$.
The ring $R$ is then Gorenstein by \ref{agree} and \ref{Goto}.
Set $M=\coprod_\bsm R/\bsm$, where the sum runs over all
maximal ideals of $R$.  For every maximal ideal 
$\bsm$ of $R$ we then obtain
$$
\dim R_\bsm=\Gpd{R_\bsm}{\big(R_\bsm/\bsm R_\bsm\big)}
\le\Gpd R{(R/\bsm)}\le\Gpd RM
$$
where the first two relations are obtained by using \ref{agree} with
\ref{ABformula}, and the third one comes from \ref{direct}.  As a
consequence, we get $\dim R\le\Gpd RM<\infty$.

(2). If $\Gpd RM<\infty$ for every finite $R$-module,
then \ref{agree} and \ref{Goto} imply that $R$ is Gorenstein, as above.

Assume now that $R$ is Gorenstein and $\bM$ is a complex
with $\HH{\!}{(\bM)}$
bounded and $\HH i{(\bM)}$ finitely generated for each $i\in\BZ$.
By \ref{semi_res} there exists a semiprojective resolution $\bP\to\bM$,
with $P_i$ finitely generated projective for each $i$.  If $\sup
\HH{\!}{(\bM)}=t$, then \ref{Goto} and \ref{agree} yield $\Gpd
R{\Ch{t}{\bP}}=n<\infty$.  Setting $g=t+n$, from the exact sequence in
the proof of Theorem \ref{Gpd_Gorenstein} we obtain $\Gpd R{\bM}\le g$.
 \end{proof}


\section{Tate cohomology}
\label{Tate cohomology}

Tate cohomology for modules over finite groups was
originally introduced, through complete resolutions, by Tate, see
\cite[Ch. XII]{ce}. It was extended to groups
of virtually finite cohomological dimension by Farrell \cite{ft} and
over groups admitting complete cohomological functors by Gedrich and
Gruenberg \cite{gg}. In various contexts cohomology theories based on
complete resolutions are constructed by Buchweitz \cite{B},
Cornick and Kropholler \cite{ck}, Avramov and Martsinkovsky.

In this section, we define and study a Tate cohomology  theory for
complexes of left modules over associative rings. Our definition is
modeled on the one of Avramov and Martsinkovsky \cite {am}.

We let $\CC_\CGp$ denote the class of complexes of finite
Gorenstein projective dimension, and we fix a complete resolution
for each complex in $\CC_\CGp$.  

Throughout this section $\bM$ denotes a complex in $\CC_\CGp$,
$\bT\xrightarrow{\tau}\bP\xrightarrow{\pi}\bM$ its chosen resolution,
and let $\bN$ denote an arbitrary complex.
For each $n\in\BZ$, the $n^{\text {th}}$ 
\emph{Tate cohomology} group is defined by
$$
\Extt{R}{n}{\bM}{\bN}=\Hh{n}{(\Hom{R}{\bT}{\bN})}.
$$
The morphism
$$
\Hom{R}{\tau}{\bN}\col\Hom{R}{\bP}{\bN}\to\Hom{R}{\bT}{\bN}
$$ 
induces for every $n\in\BZ$ a homomorphism of abelian groups
$$
\epst{R}{n}{\bM}{\bN}\col\Ext{R}{n}{\bM}{\bN}\to\Extt{R}{n}{\bM}{\bN}.
$$

\begin{5propdef}
\label{def_Tate}
\begin{enumerate}[\rm(1)]
\item  
The assignment $(\bM,\bN)\mapsto\Extt{R}{n}{\bM}{\bN}$ defines a
functor
$$
\Extt{R}{n}{-}{-}\col(\CC_\CGp)^{op}\times\CC\to\CM(\BZ)
$$ 
which is independent of choices of resolutions and liftings.
\item 
The maps $\epst{R}{n}{\bM}{\bN}$ yield a morphism of functors
$$
\epst{R}{n}{-}{-}\col\Ext{R}{n}{-}{-}\to\Extt{R}{n}{-}{-}
$$
which is independent of choices of resolutions and liftings.
\end{enumerate}
\end{5propdef}

\begin{proof}
The naturality of $\Extt{R}{n}{-}{-}$ and of $\epst{R}{n}{-}{-}$ follows
from the first part of  Proposition \ref{Tate_lifting}, applied to the
chosen complete resolutions of $\bM$ and $\bM'$. Their independence from
the choices of resolutions and liftings follows from the last part of
the same proposition.
 \end{proof}

The results below follow from the definition of Tate cohomology.

\begin{5chunk}
\label{tate_map}
If $M$ and $N$ are modules with $\Gpd R M=g<\infty$, then the natural map $\epst{R}{}{M}{N}$ is equal to $0$ for $n<0$ and is bijective for $n>g$.
\end{5chunk}

\begin{5chunk}
\label{tate_sum}
For any finite set of complexes  $\{\bM^i\}_{i\in I}$ of finite Gorenstein
projective dimension and any family of complexes $\{\bN^j\}_{j\in J}$,
there is a natural isomorphism
$$
\Extt{R}{n}{\coprod_{i\in I}\bM^i}{\prod_{j\in J}\bN^j}\cong
\prod_{(i,j)\in I\times J}\Extt{R}{n}{\bM^i}{\bN^j}\quad\text{for each}\quad n\in\BZ.
$$
\end{5chunk}
  
\begin{5chunk}
\label{tate_shift} 
There is an isomorphism
$$
\Extt{R}{n}{\Sigma^j\bM}{-}\cong\Extt{R}{n-j}{\bM}{-}\quad\text{for all}\quad n,j\in\BZ.
$$
\end{5chunk}

Unlike ordinary Ext functors, Tate cohomology is rigid. In the case
when $R$ is noetherian and $M$ is a finite module, the theorem below
specializes to \cite[(5.9)]{am}.

\begin{5theorem} 
\label{tate_pd}
Let $\bM$ be a complex with $\Gpd R\bM<\infty$. The following properties
are equivalent.
\begin{enumerate}[\rm\quad(i)]
\item $\pd R\bM<\infty$.
\item
$\Extt{R}{i}{\bM}{-}=0$ for some $i\in\BZ$.
\item[{\quad\rm(ii$'$)}]
$\Extt{R}{i}{\bM}{-}=0$ for all $i\in\BZ$.
\end{enumerate}
When $\inf\bM>-\infty$,  the properties above are also  equivalent to the following.  
\begin{enumerate}[\rm\quad(i)]
\item[\rm(iii)]
$\Extt{R}{i}{-}{\bM}=0$ for some $i\in\BZ$.
\item[{\quad\rm(iii$'$)}]
$\Extt{R}{i}{-}{\bM}=0$ for all $i\in\BZ$.
\item[\rm(iv)]
$\Extt{R}{0}{\bM}{\bM}=0$.
\end{enumerate}
\end{5theorem}

\begin{proof}  
(i)$\implies$(ii$'$). If $\pd R\bM<\infty$, then $\boldsymbol
0\xra{\tau}\bP\xra{\pi}\bM$ is a complete resolution for any
semiprojective resolution $\bP\xra{\pi}\bM$, so the group
$\Extt{R}{n}{\bM}{-}=\Hh{n}{(\Hom{R}{\boldsymbol 0}{-})}$ vanishes for all $n\in\BZ$.

(ii$'$)$\implies$(ii) is clear.

(ii)$\implies$(i). Choose a complete resolution  $\bT\xra{\tau}\bP\xra{\pi}\bM$.
Set $G=\Ch{i}{\bT}$ and let $\iota\colon G\to T_{i-1}$ be the canonical 
injection.   Since $\Extt{R}{i}{\bM}{-}=0$, the map
$\alpha=\Hom R{\iota}G$ in the commutative diagram
\[
\xymatrixrowsep{1.5pc}
\xymatrixcolsep{1.5pc}
\xymatrix{
&\Hom{R}{T_{i-1}}{G}\ar@{->}[rr]\ar@{->}[rd]_{\alpha}
&
&\Hom{R}{T_i}{G}\ar@{->}[r]
&\Hom{R}{T_{i+1}}{G}
\\
&
&\Hom{R}{G}{G}\ar@{->}[ur]
\\
&0\ar@{->}[ur]
}
\] 
is surjective.  This means that $\iota$ splits, so $G$ is projective.
Induction on $j$ shows that $\Ch{j}{\bT}$ is projective for all $j\ge
i$.  By definition, $\Ch{j}{\bP}\cong\Ch{j}{\bT}$  for all  $j\ge\Gpd
R\bM$, so $\Ch{j}{\bP}$ is projective for all $j\ge\maxx\{i,\Gpd R\bM\}$.
 From \ref{projdim} we get $\pd R\bM<\infty$.

For the rest of the proof we assume  $\inf\bM>-\infty$. 
 
(i)$\implies$(iii$'$).   By \ref{semi_res}  and \ref{projdim} choose  a  bounded semiprojective resolution $\bP\to\bM$.  Now, apply \ref{proj_exact} and \ref{total_vanishing}. 

(iii$'$)$\implies$(iii) is clear.

(iii)$\implies$(iv) follows from the isomorphism
$\Extt{R}{i}{\Sigma^{i}\bM}{\bM}\cong\Extt{R}{0}{\bM}{\bM}$.

(iv)$\implies$(i).  By \ref{semi_res} choose a bonded below semiprojective resolution $\bP\to\bM$ and  by Theorem \ref{complete} choose a complete resolution  $\bT\xra{\tau}\bP\xra{\pi}\bM$ with  $\tau_j=\id{T_j}$ for all $j\gg 0$.
By \ref{proj_exact} we get
$\Extt{R}{0}{\bM}{\bM}\cong\Hh{0}{\Hom{R}{\bT}{\bP}}=0$.
As $\tau$ is in $\Zh{0}{(\Hom{R}{\bT}{\bP})}$, there exists a
$\sigma\in\Hom{R}{\bT}{\bP}_1$ such that $\partial(\sigma)=\tau$ in
$\Hom{R}{\bT}{\bP}$. Since   $\tau_j=\id{T_j}$ for all $j\gg
0$, we then have
$$
\sigma_{j-1}\partial^\bT_j+\partial^\bP_{j+1}\sigma_j=
\id{\bT_j}\quad\text{for all}\quad j\gg 0.
$$
Since the complex $\bT$ is exact, we obtain in particular, 
$$
\partial^\bT_{j+1}\sigma_j(x)=x\quad\text{for all}\quad
x\in\im{(\partial^\bT_{j+1})}\quad\text{and all}\quad j\gg 0.
$$
Thus, the map $T_{j+1}\to\im{(\partial^\bT_{j+1})}$ splits
for all $j\gg 0$, so $\pd R\bM<\infty$ by \ref{projdim}.
\end{proof}

There exist long exact sequences of Tate cohomology groups associated
to short exact sequences of complexes in either argument.

The proof of the next result parallels that of \cite[(5.4)]{am}
and is omitted.

\begin{5proposition}
\label{tate_covariant}
For each  complex $\bM$ with  $\Gpd{R}{\bM}<\infty$ and  each exact sequence 
\[
\bN_\bu=0\to \bN\xrightarrow{\vartheta}\bN'\xrightarrow{\vartheta'}\bN''\to 0
\]
of complexes of $R$-modules there exist natural in $\bM$ and $\bN_\bu$
homomorphisms $\ett{n}{R}{\bM}{\bN_\bu}$, such that the sequence below is exact
\[
\xymatrixrowsep{0.4pc}
\xymatrixcolsep{3.7pc}
\xymatrix{
\cdots\ar@{->}[r]
&\Extt{R}{n}{\bM}{\bN}\ar@{->}[r]^-{\Extt{R}{n}{\bM}{\vartheta}}
&\Extt{R}{n}{\bM}{\bN'}\ar@{->}[r]^-{\Extt{R}{n}{\bM}{\vartheta'}}
&\Extt{R}{n}{\bM}{\bN''}
\\
{\hphantom{\cdots}}\ar@{->}[r]^-{\ett{n}{R}{\bM}{\bN_\bu}}
&\Extt{R}{n+1}{\bM}{\bN}\ar@{->}[r]
&\hskip1.9pc\cdots\hskip1.9pc
}
\]
and for each $n\in\BZ$ there is an equality
\begin{xxalignat}{3}
&\hphantom{\square}
&\ett{n}{R}{\bM}{\bN_\bu}\circ\epst{R}{n}{\bM}{\bN''}
&=\epst{R}{n+1}{\bM}{\bN}\circ\etH{n}{R}{\bM}{\bN_\bu}
&&\square
\end{xxalignat}
\end{5proposition}

To establish the existence of a long exact sequence in the first
argument, one needs a Horseshoe Lemma result for complete resolutions.
In the case of modules, such a result is proved in \cite[(5.5)]{am}.
Different arguments are needed to establish it for complexes, so a
complete proof is given.

\begin{5proposition}
\label{Tate_horseshoe}
If $0\to \bM\xra{\mu}\bM'\xra{\mu'}\bM''\to 0$ is an exact sequence of
complexes of finite Gorenstein projective dimension, then 
there exists a commutative diagram with exact rows
\begin{equation}
\tag{$*$}
\begin{gathered}
\xymatrixrowsep{2pc}
\xymatrixcolsep{2pc}
\xymatrix{
0\ar@{->}[r]
&\bM\ar@{->}[r]^{\mu}
&\bM'\ar@{->}[r]^{\mu'}
&\bM''\ar@{->}[r]
&0
\\
0\ar@{->}[r]
&\bP\ar@{->}[r]^{\ov\mu}\ar@{->}[u]^{\pi}
&\bP'\ar@{->}[r]^{\ov\mu'}\ar@{->}[u]^{\pi'}
&\bP''\ar@{->}[r]\ar@{->}[u]^{\pi''}
&0
\\
0\ar@{->}[r]
&\bT\ar@{->}[r]^{\wh\mu}\ar@{->}[u]^{\tau}
&\bT'\ar@{->}[r]^{\wh\mu'}\ar@{->}[u]^{\tau'}
&\bT''\ar@{->}[r]\ar@{->}[u]^{\tau''}
&0
}
\end{gathered}
\end{equation}
whose columns are surjective complete resolutions.
\end{5proposition}

\begin{proof}
Set $g=\max(\Gpd{R}{\bM},\Gpd{R}{\bM''})$.  Proposition \ref{semiproj_exact}
and Theorem \ref{complete} provide all  complexes and maps from the diagram $(*)$, except the complex $\bT'$ and the homomorphisms into and out of it. By Theorem \ref{complete}, we can suppose that
$$\tau_i=\id{T_i}\ \text{and}\ \tau''_i=\id{T''_i}\qquad\text{for all}\ i\geq g.$$

The rest of the proof  proceeds in
two  steps.  

\begin{step}
There exists an exact sequence $0\to\bT\to\bT'\to\bT''\to0$ of complexes
of $R$-modules, where $\bT'$ is totally acyclic with $\bT'_{\geq g}=\bP'_{\geq g}$.
 \end{step}

We start the proof of this assertion by remarking that any complex $\bT'$
appearing in an exact sequence as above is necessarily totally acyclic.
Indeed, acyclicity and the projectivity of each $T_i$ are clear.
By \ref{resolving}, for each $i\in\BZ$ the exact sequence
$$
0\to\Ch{i}{\bT}\to\Ch{i}{\bT'}\to\Ch{i}{\bT''}\to 0
$$
implies $\Ch{i}{\bT'}$ is Gorenstein projective, so
$\bT'$ is totally acyclic by Lemma \ref{syz_exact}.

Next, we describe the construction of the exact sequence in Step 1.
Theorem \ref{Gpd_ineq} gives $\Gpd{R}{\bM'}\leq g$, so we obtain
$$
\HH{i}{(\bP)}=\HH{i}{(\bP')}=\HH{i}{(\bP'')}=0\qquad \text{for all}\ i>g.
$$
Using Theorem \ref{complete} we get an exact sequence of Gorenstein
projective modules
$$
0\to\Ch{i}{\bP}\to\Ch{i}{\bP'}\to\Ch{i}{\bP''}\to 0\qquad \text{for all}\ i>g.
$$
Set $T'_i=P'_i$ for all $i\ge g$ and $\dd^{T'}_i=\dd^{P'}_i$ for all $i>g$,
and let $j$ be an integer such that $j<g$.  By descending induction, we may
assume that a commutative diagram
\[
\xymatrixrowsep{2pc}
\xymatrixcolsep{2pc}
\xymatrix{
0\ar@{->}[r]
&T_{j+1}\ar@{->}[r]\ar@{->}[d]_{\dd^{\bT}_{j+1}}
&T'_{j+1}\ar@{->}[r]\ar@{->}[d]_{\dd^{\bT'}_{j+1}}
&T''_{j+1}\ar@{->}[r]\ar@{->}[d]_{\dd^{\bT''}_{j+1}}
&0
\\
0\ar@{->}[r]
&T_j\ar@{->}[r]^{}\ar@{->}[d]_{\varepsilon_j}
&T'_j\ar@{->}[r]^{}\ar@{->}[d]_{\varepsilon'_j}
&T''_j\ar@{->}[r]\ar@{->}[d]_{\varepsilon''_j}
&0
\\
0\ar@{->}[r]
&\Ch{j}{\bT}\ar@{->}[r]^{}\ar@{->}[d]
&W_j\ar@{->}[r]^{}\ar@{->}[d]
&\Ch{j}{\bT''}\ar@{->}[r]^{}\ar@{->}[d]
&0
\\
&0
&0
&0
&
}
\]
with exact rows and columns
has been constructed.  By \ref{Gp_Ext}, the functor $\Hom{R}{-}{T_{j-1}}$
transforms  the bottom exact sequence into a short exact sequence
$$
0\to\Hom{R}{\Ch{j}{\bT''}}{T_{j-1}}\xra{}\Hom{R}{W_j}{T_{j-1}}\xra{}\Hom{R}{\Ch{j}{\bT}}{T_{j-1}}\to 0.
$$
Therefore, there exists a map $\beta\col W_j\to T_{j-1}$ such that
$\beta\alpha=\iota_{j}$, where $\iota_{j}\col\Ch{j}{\bT}\to T_{j-1}$
is the canonical injection, see the diagram below. Setting
$$
T'_{j-1}=T_{j-1}\oplus T_{j-1}''
\qquad\text{and}\qquad
\iota_{j}'(x)=
(\beta(x),\iota_{j}''\alpha'(x))\quad\text{for all}\quad x\in W_{j}
$$ 
we obtain a commutative diagram

\begin{equation}
\begin{gathered}
\tag{$**$}
\xymatrixrowsep{2pc}
\xymatrixcolsep{2pc}
\xymatrix{
&
0\ar@{->}[d]
&
0\ar@{->}[d]
&
0\ar@{->}[d]
&
\\
0\ar@{->}[r]
&\Ch{j}{\bT}\ar@{->}[r]^{\alpha}\ar@{->}[d]_{\iota_{j}}
&W_j\ar@{->}[r]^{\alpha'}\ar@{->}[d]_{\iota_{j}'}\ar@{.>}[dl]_{\beta}
&\Ch{j}{\bT''}\ar@{->}[r]\ar@{->}[d]_{\iota_{j}''}
&0
\\
0\ar@{->}[r]
&T_{j-1}\ar@{->}[r]^{}\ar@{->}[d]_{\varepsilon_{j-1}}
&T'_{j-1}\ar@{->}[r]^{}\ar@{->}[d]_{\varepsilon'_{j-1}}
&T''_{j-1}\ar@{->}[r]\ar@{->}[d]_{\varepsilon''_{j-1}}
&0
\\
0\ar@{->}[r]
&\Ch{j-1}{\bT}\ar@{->}[r]\ar@{->}[d]
&W_{j-1}\ar@{->}[r]\ar@{->}[d]
&\Ch{j-1}{\bT''}\ar@{->}[r]\ar@{->}[d]
&0
\\
&0
&0
&0
&
}
\end{gathered}
\end{equation}
where $W_{j-1}=\Coker\iota'_j$ and $\varepsilon'_{j-1}$ is the canonical
projection.  The columns and the top two rows are exact by construction,
hence so is the third row, by the Snake Lemma.  To complete the step of
the induction, set $\dd^{\bT'}_{j}=\varepsilon'_{j}\iota_{j}'$.

\begin{step}
There exists a surjective morphism $\tau'$, such that $\tau'_i=\id{T'}_i$ for 
all $i\ge g$ and the lower part of diagram $(*)$ commutes.
\end{step}

Note first that if $(*)$ commutes, then $\tau'$ is surjective because $\tau$ and
$\tau''$ are.

As the modules $P_i$ and $T''_i$ are projective, we may choose $R$-linear
splitting
$$
T_i'=T_i\oplus T_i''
\quad\text{and}\quad
P_i'=P_i\oplus P_i''
\quad\text{for all}\quad i\in\BZ
$$
such that the maps $\ov\mu_i$ and $\wh\mu_i$ respectively,
$\ov\mu'_i$ and $\wh\mu'_i$, are the canonical injections and respectively, the canonical projections. By
\cite[(4.6), p.18]{iv} there are then morphisms of complexes
$$
\theta^\bP\col\bP''\to\Sigma\bP
\qquad\text{and}\qquad
\theta^\bT\col\bT''\to\Sigma\bT
$$
such that the differentials of $\bP'$ and $\bT'$ are given, respectively, by  
$$
\partial^{\bP'}=\left(\begin{array}{cc}\partial^\bP &\theta^\bP\\
                                                                                  0& \partial^{\bP''}
                                   \end{array}\right)
\qquad\text{and}\qquad 
\partial^{\bT'}=\left(\begin{array}{cc}\partial^\bT &\theta^\bT\\
                                                                                  0 &\partial^{\bT''}\end{array}\right)
$$                                                                                
Assume that we have found a homomorphism $\Delta\colon T''\to P$,
satisfying
\begin{equation}
\tag{$1.i$}
\Delta_{i-1}\partial^{\bT''}_{i}-\partial^\bP_{i}\Delta_{i}=
\theta^\bP_i\tau''_i-\tau_{i-1}\theta^\bT_i
\end{equation}
for all $i\in\BZ$.
As in the classical case in \cite[(V.2.3)]{ce},
a simple computation shows that the map $\tau'\col\bT'\to\bP'$ defined by 
$$
\tau'=\left(\begin{array}{cc}\tau&\Delta\\
                                                   0&\tau'' \end{array}\right)
                                                   $$
is a morphism making the diagram $(*)$  commutative.

We produce $\Delta_i$ by descending induction on $i$.  As $\tau_i$
and $\tau''_i$ are identity maps for $i\ge g$, for these $i$ 
we set $\Delta_i=0$. Note that $\tau'_i=\id{T'_i}$ for all $i\geq g$.
Let $j$ be an integer such that $j<g$ and assume, by induction, that
$\Delta_j$ has been constructed and (1.$j+1$) holds.  We then have
\begin{multline*}
(\theta^\bP_j\tau''_j-\tau_{j-1}\theta^\bT_j+\partial^\bP_{j}\Delta_{j})\partial_{j+1}^{\bT''}=\\
{\begin{aligned}
&{=\theta^\bP_j\tau''_j\partial_{j+1}^{\bT''}-\tau_{j-1}\theta^\bT_j\partial_{j+1}^{\bT''}
+\partial^\bP_j(\theta^\bP_{j+1}\tau''_{j+1}-\tau_{j}\theta^\bT_{j+1}+\partial^\bP_{j+1}\Delta_{j+1})}\\
&{=(\theta^\bP_j\tau''_j\partial_{j+1}^{\bT''}+\partial^\bP_j\theta^\bP_{j+1}\tau''_{j+1})
-(\tau_{j-1}\theta^\bT_j\partial_{j+1}^{\bT''}+\partial^\bP_j\tau_{j}\theta^\bT_{j+1})}\\
&{=(\theta^\bP_j\partial^{\bP''}_{j+1}+\partial^\bP_j\theta^\bP_{j+1})\tau''_{j+1}-
\tau_{j-1}(\theta^\bT_j\partial_{j+1}^{\bT''}+\partial^\bT_j\theta^\bT_{j+1})}\\
&{=0.}
\end{aligned}}
\end{multline*}
Since the complex $\bT''$ is totally acyclic, for $P=P_{j-1}$ we get an
exact sequence
$$
{\Hom{R}{T''_{j-1}}{P}
\xra{\Hom{R}{\partial^{\bT''}_{j}}{P}}\Hom{R}{T''_j}{P}
\xra{\Hom{R}{\partial^{\bT''}_{j+1}}{P}}\Hom{R}{T''_{j+1}}{P}}
$$
Thus, there exists a homomorphism $\Delta_{j-1}\col T''_{j-1}\to P_{j-1}$ satisfying
\[
\Delta_{j-1}\partial_{j}^{\bT''}=
\theta^\bP_j\tau''_j-\tau_{j-1}\theta^\bT_j+\partial^\bP_{j}\Delta_{j}
\]
that is, formula (1.$j$) is satisfied. This finishes the step of the induction.
\end{proof}

Now, the proof of the following proposition is parallel to that of \cite[(5.6)]{am}.

\begin{5proposition}
\label{tate_contravariant}
 For each exact sequence  
\[
\bM_\bu=0\to \bM\xrightarrow{\mu}\bM'\xrightarrow{\mu'} \bM''\to 0
\]
of complexes of finite Gorenstein projective dimension and for every complex of $R$-modules $\bN$
there exist natural in $\bM_\bu$ and $\bN$ homomorphisms
$\ett{n}{R}{\bM_\bu}{\bN}$, such that the sequence below is exact
\[
\xymatrixrowsep{0.4pc}
\xymatrixcolsep{3.7pc}
\xymatrix{
\cdots\ar@{->}[r]
&\Extt{R}{n}{\bM''}{\bN}\ar@{->}[r]^{\Extt{R}{n}{\mu'}{\bN}}
&\Extt{R}{n}{\bM'}{\bN}\ar@{->}[r]^{\Extt{R}{n}{\mu}{\bN}}
&\Extt{R}{n}{\bM}{\bN}
\\
{\hphantom{\cdots}}\ar@{->}[r]^{\hspace{-1cm}\ett{n}{R}{\bM_\bu}{\bN}}
&\Extt{R}{n+1}{\bM''}{\bN}\ar@{->}[r]
&\hskip1.9pc\cdots\hskip1.9pc
}
\]
and for each $n\in\BZ$ there is an equality
\begin{xxalignat}{3}
&\hphantom{\square}
&\ett{n}{R}{\bM_\bu}{\bN}\circ\epst{R}{n}{\bM}{\bN}
&=\epst{R}{n+1}{\bM''}{\bN}\circ\etH{n}{R}{\bM_\bu}{\bN}
&&\square
\end{xxalignat}
\end{5proposition}

\section{Complexes with bounded below homology}
\label{Complexes with bounded below homology}

Gorenstein dimension of complexes with bounded below homology admits
additional descriptions.  In this section we concentrate on such
complexes whose Gorenstein projective dimension is finite.  First, we
note that they form a proper subclass of the class of complexes of finite
Gorenstein projective dimension.

\begin{5example}  
The complex $\bM=\coprod_{i\leq 0}(\Sigma^i R)$
satisfies 
$$
\inf\HH{\!}{(\bM)}=-\infty\qquad\text{and}\qquad\Gpd R\bM=0.
$$

Indeed, the first equality is clear. For the second, 
apply Corollary \ref{coproducts} and \ref{shiftt}.
\end{5example}

\begin{5chunk}
A \emph{Gorenstein projective resolution}
of $\bM$ is  a complex of Gorenstein projective modules $\bG$  such
that $\bG\simeq\bM$. Such a resolution is \emph{finite} if $G_i=0$
for all $|i|\gg 0$; it  is  \emph{special} if it is finite, $\inf\bG=\inf\HH{\!}{(\bM)}$, and $G_i$
is projective for all $i>\inf\HH{\!}{(\bM)}$.
\end{5chunk}

The notion of special Gorenstein projective resolution is an extension
to complexes of the notion of finite strict resolution which was defined
in the case of finite modules over commutative rings by Avramov and
Martsinkovsky in  \cite{am}, see also Remark \ref{strict-special}. It
is used in  Section \ref{Relative} to define  relative cohomology
for arbitrary modules over associative rings.

When given a complex $\bM$ with bounded below homology, of
finite Gorenstein projective dimension, one can easily construct
a special Gorenstein projective resolution, see Construction
\ref{Gpd_special_constr}. The next result, which is proved at the end
of this section, shows that in the case of bounded below complexes one
can compute the Gorenstein  projective dimension using only special
resolutions.

\begin{5theorem}
\label{Gpd_finite_special}
If $\bM$ is a complex with $\inf\HH{\!}{(\bM)}>-\infty$, then there are equalities 
\begin{align*}
\Gpd{R}{\bM}
=&\inf\left\{\sup\bG\left|\begin{gathered}
                         \bG\ \text{is a finite Gorenstein}\\ 
                          \text{projective resolution of}\ \bM         
                          \end{gathered}
                     \right.
     \right\}\\
=&\inf\left\{\sup\bG\left|\begin{gathered}\bG\ \text{is a
special Gorenstein}\\ \text{ projective resolution of}\
\bM\end{gathered}\right.\right\}
\end{align*}
\end{5theorem}
 
In the next remark, we compare our notion of Gorenstein projective
dimension for complexes with earlier concepts, which are restricted
to special classes of complexes.

\begin{5remark}
\label{old_new}
When $R$ is commutative and noetherian,  Yassemi defines in \cite{y}
a Gorenstein dimension for complexes $\bM$ such that $\HH{\!}{(\bM)}$
is bounded and $\HH{i}{(\bM)}$ is finite for all $i\in\BZ$.  He says that
$\bM$ has finite  Gorenstein dimension if H({\bf R}$\Hom{R}{\bM}{R}$)
is finite and the canonical map in the derived category
$$
\bM\to\text{{\bf R}Hom}_R({\text{\bf R}}\Hom{R}{\bM}{R},R)
$$ 
is an isomorphism. When this is the case, he calls the number
$$
-\inf(\HH{\!}{(\text{\bf R}\Hom{R}{\bM}{R}))}
$$
the Gorenstein dimension of $R$.  When $M$ is a finite $R$-module,
he shows in \cite[(2.7)]{y} that it is equal to the Gorenstein dimension of
Auslander-Bridger \cite[(3.8)]{ab}.

When $\HH{\!}{(\bM)}$ is bounded below, Christensen defines in
\cite[(4.4.2)]{cl} the Gorenstein projection dimension of $\bM$  to be
$$
\inf\left\{\sup\bG\left|\begin{gathered}
                         \bG\ \text{ is a complex of}\\
\text{Gorenstein projective modules with}\\ 
                           \inf\bG>-\infty\quad   \text{and}\quad \bG\simeq\bM        
                          \end{gathered}
                     \right.
     \right\}     
$$
He shows in \cite[(2.3.8)]{cl} that if the hypotheses of Yassemi's 
definition are satisfied, then the two Gorenstein dimensions coincide.

The first equality in Theorem \ref{Gpd_finite_special} shows that for
bounded below complexes our definition of Gorenstein projective dimension
coincides with that of Christensen.
\end{5remark}

For the rest of this section, unless otherwise specified, $\bM$ is
a complex with
$$
\inf\HH{\!}{(\bM)}=t<\infty.
$$

\begin{5construction}
\label{Gpd_special_constr} 
Suppose $\Gpd R\bM=g<\infty$. By \ref{semi_res} and Theorem \ref{complete} choose a surjective complete resolution $\bT\xra{\tau}\bP\xra{\pi}\bM$ with $\inf\bP=t$
and $\tau_i$ bijective for all $i\ge g$. Set $\bL=\Ker\tau$ and form a
complex $(\bG,\partial^{\bG})$ by setting
\begin{alignat*}{3}
G_i = \begin{cases} 0&\text{for}\ i<t\,;\\
                    \Ch{t}{\bT}&\text{for}\  i=t\,;\\
                    L_{i-1} &\text{for}\ i>t\,; 
      \end{cases}
\qquad
&\text{and}&
\qquad
\partial_i^{\bG}=\begin{cases}0&\text{for}\ i<t+1\,;\\
                              -\varsigma&\text{for}\ i=t+1\,;\\ 
                              -\partial_{i-1}^{\bL}&\text{for}\ i>t+1\,,
                 \end{cases}
\end{alignat*}
where $\varsigma\col L_t\to\Ch{t}{\bT}$ is the natural homomorphism of modules. 

Let $\iota\col\bG\to\Sigma\bL$ be the morphism of complexes given by
$$
\iota_i=\begin{cases}
              0&\text{for}\ i<t\,;\\
              \varsigma'&\text{for}\ i=t\,;\\
              \id {L_{i-1}}&\text{for}\ i>t\,,        
              \end{cases}
$$
where $\varsigma'\col \Ch{t}{\bT}\to T_{t-1}=L_{t-1}$ is the canonical inclusion.

Let $\lambda\col\bP\to\Sigma\bL$  be the morphism defined in
\ref{ker_complete} via \ref{theta_construction} and  define a map $\gamma\col\bP\to\bG$  by the
formulas
$$
\gamma_i=\begin{cases}0&\text{for}\ i<t\,;\\
                              \varsigma''&\text{for}\ i=t\,;\\ 
                              \lambda_i&\text{for}\ i>t\,,
                 \end{cases}
$$
where $\varsigma''\col P_t\to\Ch{t}{\bT}$ is the natural homomorphism of
modules induced by $\lambda_t$.
 \end{5construction}

\begin{5lemma}  
\label{Gpd_special_lemma}
With the notation above, the  diagram  
\[
\xymatrixrowsep{2pc}
\xymatrixcolsep{2pc}
\xymatrix{
\bG\ar@{->}[d]_{\iota}
&\bP\ar@{->}[l]_{\gamma}\ar@{=}[d]\ar@{->}[r]^{\pi}
&\bM\ar@{=}[d]
\\
\Sigma\bL
&\bP\ar@{->}[l]_{\lambda}\ar@{->}[r]^{\pi}
&\bM
}
\]
is commutative and all arrows are quasiisomorphisms.  In particular, $\bG$
is a special Gorenstein projective resolution of $\bM$ with $\sup\bG=g$.

When $\bM$ is a module, $M$, there is a quasiisomorphism
$\varepsilon\col\bG\to M$ with $\varepsilon\gamma=\pi.$
 \end{5lemma}

\begin{proof} 
Commutativity follows from the definitions of the maps.  The definition
of $\bG$, the equality $L_{\leq t-1}=T_{\leq t-1}$ and the exactness of
$\bT$ yield  
$$
\HH{i}{(\bG)}=\begin{cases} 0&\text{for}\ i<t\,;\\
                           \HH{i-1}{(\bL)}&\text{for}\ i\geq t\,.
               \end{cases}            
$$           
Therefore, $\iota$ is a quasiisomorphism.  The morphism $\lambda$
is a quasiisomorphism by \ref{ker_complete}, hence from the equality
$\iota\gamma=\lambda$  we obtain that $\gamma$ is a quasiiomorphism. We
have $\sup\bG=g$ by Construction \ref{Gpd_special_constr}.  

If $\bM=M$ is a module, then $t=0$,  $\bP$ becomes a projective resolution
of $M=\Ch{0}{\bP}$ and $\tau_0\col T_0\to P_0$ induces a map
$\ov\tau_0\col\Ch{0}{\bT}\to\Ch{0}{\bP}$.  Setting
$\varepsilon_i=0$ for $i\not=0$ and $\varepsilon_0=\ov\tau_0$
we get a quasiisomorphism $\varepsilon\col\bG\to M$ with
$\varepsilon\gamma=\pi$.
\end{proof}

\begin{5construction}
\label{special_Gpd_constr} 
Suppose that $\bG\simeq\bM$ is a special Gorenstein projective
resolution.  By definition, the module $G_t$ is  Gorenstein
projective, so there exists a totally acyclic complex $\bS$ with
$\Ch{t+1}{\bS}=G_t$. 

Let $(\bK,\partial^{\bK})$ denote the complex of
projectives with
\begin{alignat*}{3}
K_i=\begin{cases} S_i&\text{for}\ i\leq t\,;\\
                      G_i&\text{for}\ i> t\,;
         \end{cases}
\qquad
&
\text{and}
&
\qquad
\partial^{\bK}_i=\begin{cases}\partial^\bS_i&\text{for}\ i<t+1\,; \\
                                 \varepsilon\circ\partial^\bG_{t+1}&\text{for}\ 
                                               i=t+1\,;\\
                                 \partial^\bG_i&\text{for}\
                                               i>t+1\,,
                    \end{cases}
\end{alignat*}    
where $\varepsilon\col G_t\to S_t$ is the canonical injection.

Let $\iota'\col\bG\to\bK$ be the morphism of complexes given by
$$
\iota'_i=\begin{cases}
              0&\text{for}\ i<t\,;\\
              \varepsilon&\text{for}\ i=t\,;\\
              \id {G_{i}}&\text{for}\ i>t\,.        
              \end{cases}
$$
\end{5construction}

\begin{5lemma} 
\label{special_Gpd_lemma}
With the notation above, for any semiprojective resolution
$\bP\xra{\pi}\bM$, there exists a commutative diagram of quasiisomorphisms
\[
\xymatrixrowsep{2pc}
\xymatrixcolsep{2pc}
\xymatrix{
\bG\ar@{->}[d]_{\iota'}
&\bP\ar@{->}[l]_{\gamma}\ar@{=}[d]\ar@{->}[r]^{\pi}
&\bM\ar@{=}[d]
\\
{\bK}
&\bP\ar@{->}[l]_{\kappa}\ar@{->}[r]^{\pi}
&\bM
}
\]
If $\inf\bP=t<\infty$, then there exists a surjective complete resolution
$$
\bT\xra{\tau}\bP\xra{\pi}\bM
$$
with $\Ker\tau=\Sigma^{-1}\bK$.
\end{5lemma}

\begin{proof} 
A quasiisomorphism $\gamma$ is obtained from \ref{semiprojective_quism}. It
is easy to see that $\iota'$ is a quasiisomorphism.  To complete the
square on the right hand side, set $\kappa=\iota'\gamma$.

Assume $\inf\bP=t<\infty$.  Set $\bT=\Sigma^{-1}\Cone(\kappa)$; this
is a complex of projectives, and it is exact because $\kappa$ is a
quasiisomorphism.  For all $i\le t$ we have $\Ch{i}{\bT}=\Ch{i}{\bS}$
by construction, so the module $\Ch{i}{\bT}$ is Gorenstein projective;
Lemma \ref{syz_exact} now shows that $\bT$ is totally acyclic. The
exact sequence of complexes
$$
0\to\Sigma^{-1}\bK\to\bT\xra{\tau}\bP\to0
$$
gives a surjective complete resolution with the desired property.
\end{proof}

{\bf Remark.} The reader might note that Constructions  \ref{Gpd_special_constr} and \ref{special_Gpd_constr} are
"inverse to each other".

\begin{proof}[Proof of Theorem \emph{\ref{Gpd_finite_special}}] 
If $\inf\HH{\!}{(\bM)}=\infty$, then $\bM$ is exact and ${\bf 0}\cong\bM$ is
a finite (special) Gorenstein projective resolution, therefore the infima
of the sets in the right hand side are $-\infty$. By \ref{Gpd_exact},
$\Gpd R\bM=-\infty$, so the theorem holds in this case. Thus, we may
assume that $\inf\HH{\!}{(\bM)}$ is finite.

Set  $\Gpd{R}{\bM}=g$ and let $h$
(respectively, $h'$) denote the number on the right hand side of the
first (respectively, second) equality of the theorem.

We show first $g\leq h$. If $h=\infty$ it is clear. If $h<\infty$, then there exists  a finite Gorenstein
projective resolution $\bG$ of $\bM$; set $s=\sup\bG$. Let $\bP\to\bM$
be a semiprojective resolution with $P_i=0$ for all $i<t$
and let $\gamma\col\bP\to\bG$ be a quasiisomorphism given by
\ref{semiprojective}. The mapping cone $\Cone(\gamma)$ is an exact
complex of Gorenstein projective modules, bounded below, so every syzygy
is a Gorenstein projective module by \ref{resolving}, in particular so is
$\Ch{s+1}{\Cone(\gamma)}=\Ch{s}{\bP}$. As $\bG$ is quasiisomorphic to
$\bM$, $\sup\HH{\!}{(\bM)}\leq s$. From Theorem \ref{complete}
we get $g\leq s$, so $g\leq h$.

The inequality $h\leq h'$ is clear.

We show next  $h'\leq g$. If $g=\infty$ it is clear. If $g<\infty$, then by  Lemma \ref{Gpd_special_lemma} there exists a special Gorenstein projective
resolution $\bG$ of $\bM$ such that $\sup\bG=g$.
\end{proof}

\section{Relative cohomology for modules}
\label{Relative} 

In this section $R$ is an associative ring and $M$ is a left $R$-module.

When $R$ is two-sided noetherian and $M$ admits a ``proper'' resolution
by finite modules of Gorenstein dimension zero, Avramov and Martinkovsky
\cite{am} associate to $M$ relative cohomology groups.  They study the
behavior of these groups with respect to short exact sequences, and relate
them to the absolute cohomology groups and to the Tate cohomology groups.

Holm \cite{h} extends the definition of relative cohomology groups
to the case when $M$ admits a proper resolution by Gorenstein projective
modules over an arbitrary ring $R$.  He focuses on the relation of 
these relative groups with those defined using a proper resolution of the
second argument by Gorenstein injective modules.

Our purpose here is to show that the result of \cite{am} extend to the
general setup of \cite{h}.  As most arguments carry over, we just
give indications of proofs.  In Remark \ref{explanation} we comment on
obstacles to defining relative cohomology groups for complexes.

\begin{5chunk}
A complex $\bC$ is called \emph{proper exact} if  $\Hom{R}{E}{\bC}$
is exact for all Gorenstein projective modules $E$.  A \emph{proper
Gorenstein projective resolution} of $M$ is a complex of Gorenstein
projective modules $\bG$ together with a morphism $\varepsilon\col \bG\to
M$, such that the  complex $\Cone(\varepsilon)$ is  proper exact.
 \end{5chunk}

Let $\CGp$ denote the class of Gorenstein projective
modules and let $\ov\CGp$ denote the class of modules that admit some
proper resolution.

\begin{5remark}
\label{strict-special}
Let $M$ be an $R$-module with $\Gpd RM<\infty$. By Lemma
\ref{Gpd_special_lemma} there exists a special Gorenstein projective
resolution $\bG\xla{\gamma}\bP\xra{\pi}M$ and a quasiisomorphism
$\varepsilon\col\bG\to M$ such that $\varepsilon\gamma=\pi$.  Such a
resolution is called finite strict in \cite{am} and any finite strict
resolution is proper by \cite[(4.1.3)]{am}. Thus, every module of finite
Gorenstein projective dimension is in $\ov\CGp$.  \end{5remark}

\begin{5chunk}
\label{relative_def}
For each $M\in\ov\CGp$ choose a proper resolution $\bG\to M$. For each $R$-module
$N$ and each $n\in\BZ$, the $n^{\text th}$ \emph{relative cohomology} group is
defined by 
$$
\Extr{n}{M}{N}=\Hh{n}{\Hom{R}{\bG}{N}}.
$$
Choose a projective resolution $\bP\to M$ and a morphism
$\gamma\col\bP\to\bG$ lifting the identity on $M$.
The morphism of complexes of abelian groups
$$
\Hom{R}{\gamma}{N}\col\Hom{R}{\bG}{N}\to\Hom{R}{\bP}{N}
$$
induces for every $n\in\BZ$ a natural homomorphism of abelian groups
$$
\epsr{n}{M}{N}\col\Extr{n}{M}{N}\to\Ext{R}{n}{M}{N}
$$ 
The groups and maps defined above do not depend on the choices of resolutions
and liftings used in their constructions, see \cite[(4.2)]{am} or \cite[(2.4)]{h}.
 \end{5chunk}

The following two propositions correspond to \cite[(4.6)]{am}, \cite[(4.4)]{am},
whose proofs extend verbatim to the more general setting below.
 
\begin{5proposition}
\label{rel_covariant}
For each proper exact sequence $\bM=0\to M\xra{\mu}M'\xra{\mu'}M''\to0$
of $R$-modules in $\ov\CGp$ and each $R$-module $N$ there
are natural in $\bM$ and $N$ homomorphisms $\etr{n}{\bM}{N}$,
such that the sequence below is exact
\[
\xymatrixrowsep{0.4pc}
\xymatrixcolsep{3.7pc}
\xymatrix{
\cdots\ar@{->}[r]
&\Extr{n}{M''}{N}\ar@{->}[r]^-{\Extr{n}{\mu'}{N}}
&\Extr{n}{M'}{N}\ar@{->}[r]^-{\Extr{n}{\mu}{N}}
&\Extr{n}{M}{N}
\\
{\hphantom{\cdots}}\ar@{->}[r]^-{\etr{n}{\bM}{N}}
&\Extr{n+1}{M''}{N}\ar@{->}[r]
&\hskip1.9pc\cdots\hskip1.9pc
}
\]
and for each $n\in\BZ$ there is an equality
\begin{xxalignat}{3}
&\hphantom{\square}
&\etH{n}{R}{\bM}{N}\circ\epsr{n}{M}{N}
&=\epsr{n+1}{M''}{N}\circ\etr{n}{\bM}{N}
&&\square
\end{xxalignat}
\end{5proposition}

\begin{5proposition}
\label{rel long exact 2}
For each $R$-module $M\in\ov\CGp$ and each proper exact
sequence $\bN=0\to N\xra{\nu}N'\xra{\nu'}N''\to 0$ of $R$-modules there
exist natural in $M$ and $\bN$ homomorphisms $\etr{n}{M}{\bN}$,
such that the sequence below is exact
\[
\xymatrixrowsep{0.4pc}
\xymatrixcolsep{3.7pc}
\xymatrix{
\cdots\ar@{->}[r]
&\Extr{n}{M}{N}\ar@{->}[r]^-{\Extr{n}{M}{\nu}}
&\Extr{n}{M}{N'}\ar@{->}[r]^-{\Extr{n}{M}{\nu'}}
&\Extr{n}{M}{N''}
\\
{\hphantom{\cdots}}\ar@{->}[r]^-{\etr{n}{M}{\bN}}
&\Extr{n+1}{M}{N}\ar@{->}[r]
&\hskip1.9pc\cdots\hskip1.9pc
}
\]
and for each $n\in\BZ$ there is an equality
\begin{xxalignat}{3}
&\hphantom{\square}
&\etH{n}{R}{M}{\bN}\circ\epsr{n}{M}{N''}
&=\epsr{n+1}{M}{N}\circ\etr{n}{M}{\bN}
&&\square
\end{xxalignat}
\end{5proposition}

When $\Gpd RM=g<\infty$ Remark \ref{tate_map} describes the map
$\epst{R}{}{M}{N}$ for $n<0$ and $n>g$. The next result analyzes
this map for $1\leq n\leq g$.  It extends \cite[(7.1)]{am}, whose proof
also carries over.  We indicate how the original argument can be
shortened.

\begin{5theorem}
\label{relative_Tate}
Let $M$ be a module with $\Gpd RM=g<\infty$.  For each  $R$-module $N$
there exist homomorphisms $\deltat{R}{n}{M}{N}$, natural in $M$ and $N$,
such that the following sequence is exact
 \[
\xymatrixrowsep{.4pc}
\xymatrixcolsep{3 pc}
\xymatrix{
0\ar@{->}[r]
&\Extr{1}{M}{N}\ar@{->}[r]^{\epsr{1}{M}{N}}
&\Ext{R}{1}{M}{N}\ar@{->}[r]
&\qquad\cdots\\
{\phantom{0}}\ar@{->}[r]
&\Extr{n}{M}{N}\ar@{->}[r]^{\epsr{n}{M}{N}}
&\Ext{R}{n}{M}{N}\ar@{->}[r]^{\epst{R}{n}{M}{N}}
&\Extt{R}{n}{M}{N}
\\
{\phantom{0}}\ar@{->}[r]^{\deltat{R}{n}{M}{N}\qquad\quad}
&\Extr{n+1}{M}{N}\ar@{->}[r]
&\qquad\cdots\qquad\ar@{->}[r]
&\Extt{R}{g}{M}{N}\ar@{->}[r]
&0
}
\]
\end{5theorem}

\begin{proof}
As in the proof of \cite[(7.1)]{am}, but using Lemma \ref{special_Gpd_lemma}  instead of \cite[(3.8)]{am}, we obtain a split exact sequence 
\begin{equation*}
\tag{*}
\qquad 0\to\Hom{R}{\bP}{N}\to\Hom{R}{\bT^\flat}{N}\to\Hom{R}{\Sigma^{-1}\bG}{N}\to 0
\end{equation*}
Its  long cohomology sequence 
$$\cdots\xra{\deltat{R}{n-1}{M}{N}}\Extr{n}{M}{N}\xra{\eth^{n-1}}\Ext{R}{n}{M}{N}
\xra{\Hh{n}{(\Hom{R}{\tau^\flat}{N}})}\Hh{n}{(\Hom{R}{\bT^\flat}{N})}$$
induces the desired long  exact sequence of the theorem, which is natural in $M$ and $N$.

\smallskip

Here we remark that the naturality of the map $\eth^{n-1}$, follows by general theory.
Indeed,  the exact sequence (*) is obtained as follows:  we consider the cone exact sequence of the map $\gamma\col\bP\to\bG$ of Lemma \ref{special_Gpd_lemma}  and we apply $\Sigma^{-1}$ to it, so we obtain the split exact sequence 
$$0\to\Sigma^{-1}\bG\to\bT^\flat \to\bP\to 0.$$  Applying  $\Hom{R}{-}{N}$, we get (*). 
By  \ref{theta_construction} and \ref{theta_arg2} we obtain the expression
$$ \eth^{n-1}=(-1)^n\epsr{n}{M}{N}\ \text{for all}\  n>0$$ 
for the connecting homomorphism $\eth$.

\end{proof}

\begin{5remark}
\label{explanation}
Let $\bM$ be a complex with $\Gpd R{\bM}<\infty$ and
$\inf\HH{}{(\bM)}>\infty$.  By Lemma \ref{Gpd_special_lemma}, choose  a special Gorenstein projective resolution $\bG\simeq\bM$.  If $\bM=M$ and $N$ are modules, then
$\Hh{n}{(\Hom{R}{\bG}{N})}$ is the relative cohomology group
$\Extr{n+1}{M}{N}$, see \ref{relative_def}.

It is tempting to use $\Hh{n}{(\Hom{R}{\bG}{\bN})}$ to define relative
cohomology groups for complexes.  However, I do not know whether
this construction has the necessary uniqueness and functoriality.
More precisely, if $\bM'$ is a complex with $\Gpd R{\bM'}<\infty$
and $\inf\HH{}{(\bM)}\ge\inf\HH{}{(\bM')}>\infty$, then for each
morphism $\mu\col\bM\to\bM'$ there exists a morphism of complexes
$\wh\mu\col\bG\to\bG'$ such that if $\mu$ is an isomorphism (in
particular, if $\mu=\id\bM$), then so is $\Hh{n}{(\Hom{R}{\wh\mu}{\bN})}$. I do not know whether this map it depends on the choice of $\wh\mu$.
 \end{5remark}

\section*{Acknowledgments}
I should like to thank my thesis adviser Luchezar Avramov for  his patient
guidance throughout the development and writing this article.  I thank
Sean Sather-Wagstaff for pointing me to Iversen's book \cite{iv}, Anders
Frankild and Henrik Holm for comments on an earlier version,
and Srikanth Iyengar for Theorem \ref{Gpd_Gorenstein_com}.
I thank the referee for a careful reading of the paper and for 
many useful suggestions.


\begin{thebibliography}{99}


\bibitem{ab} 
 Auslander, M., Bridger, M.,
\emph{ Stable module theory\/},
Mem. Amer. Math. Soc. {\bf 94} (1969).

\bibitem{av}
Avramov, L.L.,
\emph{Homological dimensions and related invariants of modules over local rings\/},
Representations of Algebras, ICRA IX (Beijing, 2000), vol. I, Beijing
Normal Univ. Press 2002, 1--39.

\bibitem{af}
Avramov, L.L., Foxby, H.-B.,
\emph{Homological dimension of unbounded complexes,\/}
J. Pure Appl. Alg. {\bf 71} (1991), 129--155.

\bibitem{afh}
Avramov, L.L., Foxby, H.-B., Halperin, S.,
\emph{Differential graded homological algebra,\/}
 Preprint, (2003).

\bibitem{avgp} 
 Avramov, L.L., Gasharov, V.N., Peeva, I.V.,
{\em Complete intersection dimension, \/}
 Publ. Math. I.H.E.S., {\bf 86} (1997), 67--114.

\bibitem{am}
Avramov, L.L., Martsinkovsky, A.,
\emph{Absolute, relative, and Tate cohomology of modules of finite
  Gorenstein dimension,\/} Proc. London Math. Soc. (3) {\bf 85}
(2002), 393--440.

\bibitem{b}
Bourbaki, N.,
\emph{Alg\`{e}bre Chaptre 10, Alg\`{e}bre Homologique,\/}
 Masson, Paris, New York,  1980.

\bibitem{bh}
Bruns, W., Herzog, J.,
\emph{Cohen Macaylay ring,\/}
(revised edition), Advances in Mathematics, Vol. 39, Cambridge
Univ. Press, Cambridge, UK, 1996.
 
\bibitem {B}
Buchweitz, R.-O.,
{\em  Maximal Cohen-Macaulay modules and Tate cohomology over
  Gorenstein rings\/},
Preprint, Univ. Hannover, 1986.  

\bibitem{ce} 
Cartan, H., Eilenberg, S.,
\emph{Homological Algebra},
Princeton University Press, Princeton, NJ, reprint 1999.

\bibitem{ck}
Cornick, J., Kropholler, P.H.,
\emph{On complete resolutions }, 
Topology Appl. {\bf 78} (1997), 235--250.

\bibitem{cl}
Christensen, L.W.,
\emph{Gorenstein dimensions},
Lecture Notes in Math. {\bf 1747}, Springer, Berlin, 2000.

\bibitem{ej1}
Enochs, E.E., Jenda, O.M.G.,
\emph{Gorenstein injective and projective modules},
Math. Z. {\bf 220}, (1995),  611--633.

\bibitem{ej2}
Enochs, E.E., Jenda, O.M.G.,
\emph{Relative Homological Algebra},
De Gruyter Exp. Math. {\bf 30}, De Gruyter, Berlin, 2000.

\bibitem{ft} 
Farrell, F.T.,
{\em An extension of Tate cohomology to a class of infinite groups\/},
J. Pure Appl. Alg. {\bf 10} (1977), 153--161.

\bibitem{gg} Gedrich, T.V., Gruenberg, K.W.,
{\em Complete cohomological functors on groups\/},
Topology Appl. {\bf 25} (1987), 203-223.

\bibitem{g}
Gerko, A.A.,
{\em On homological dimensions\/}, 
Sb. Math. {\bf 192} (2001), 1165--1179.

\bibitem{go} 
Goto, S.,
{\emph Vanishing of $\Ext{A}{i}{M}{A}$\/},
J. Math. Kyoto Univ. {\bf 22} (1982), 481--484. 

\bibitem{h} Holm, H.,
\emph{Gorenstein derived functors},
Proc. Amer. Math. Soc., {\bf 132}  (2004), 1913--1923. 


\bibitem{hh} Holm, H.,
\emph{Gorenstein homological dimensions,\/} 
J. Pure Appl. Alg. {\bf 189} (2004), 167--193.


\bibitem{iv} Iversen, B.,
\emph{Cohomology of Sheaves,\/}
Springer-Verlag, Berlin-Heidelberg, 1986.   

\bibitem{sw} Sather-Wagstaff, S.,
\emph{Complete intersection dimension for complexes,\/}
J. Pure Appl. Alg. {\bf 190} (2004), 267--290.

\bibitem{y} Yassemi, S.,
\emph{Gorenstein dimension},
Math. Scand. {\bf 77} (1995), 161--174.

\bibitem{v} Veliche, O.,
\emph{Construction of modules of finite Gorenstein dimension,\/}
J. Algebra {\bf 250} (2002), 427--449.

\bibitem{w} Weibel, C.,A.,
\emph{An introduction to homological algebra,\/}
Cambridge Studies in advanced mathematics {\bf 38},  Cambridge University Press, Cambridge, 1994.

\end{thebibliography}
\end{document}